\documentclass[journal, twocolumn]{IEEEtran}

\IEEEoverridecommandlockouts                              % This command is only needed if 

\usepackage{tikz}
\usepackage{textcomp}

\usetikzlibrary{matrix}
\usetikzlibrary{angles,quotes}

\usepackage{graphics,graphicx} % for pdf, bitmapped graphics files
\usepackage{epsfig} % for postscript graphics files
\usepackage{mathptmx} % assumes new font selection scheme installed
\usepackage{times} % assumes new font selection scheme installed
\usepackage{amsmath} % assumes amsmath package installed
\usepackage{amssymb,nth, mathtools,dsfont}  % assumes amsmath package installed
\usepackage{amsbsy}
\usepackage{epstopdf}
\usepackage[noadjust]{cite}
\usepackage{romannum}
\usepackage{xcolor}
\usepackage{subcaption}
\usepackage{multirow}

\newtheorem{theorem}{Theorem}
\newtheorem{lemma}{Lemma}
\newtheorem{definition}{Definition}
\newtheorem{remark}{Remark}
\newtheorem{conjecture}{Conjecture}
\newtheorem{proposition}{Proposition}
\newtheorem{algorithm}{Algorithm}
\newtheorem{corollary}{Corollary}

\makeatletter
\newcommand*{\rom}[1]{\expandafter\@slowromancap\romannumeral #1@}

\title{Duality bounds for discrete-time Zames-Falb multipliers}

\author{Jingfan Zhang,  Joaquin Carrasco and William P. Heath% <-this % stops a space
	\thanks{Department of Electrical $\&$ Electronic Engineering,  School of Engineering,   University of Manchester, M13 9PL, UK. \newline 
		{\tt\small jingfan.zhang@manchester.ac.uk}\newline 
		{\tt\small joaquin.carrascogomez@manchester.ac.uk} \newline 
		{\tt\small william.heath@manchester.ac.uk}   
	}
}

\begin{document}

\maketitle

\begin{abstract}
%In the first part of the paper, we revisit duality results for Zames--Falb multipliers for continuous-time and we develop the counterpart results in discrete time. The main contribution of this paper are to develop a discrete-time counterpart of the dual bound results and to show that these bounds are numerically tractable for discrete-time. Two different solutions are provided: a closed-form solution and a linear program. In contrast with the continupus case, our result show a negligible dual gap, i.e. gap between being able to find an LTI multiplier and to discard its existence, hence a full solution to the problem of finding a suitable Zames--Falb multiplier in discrete-time is finally concluded. 
We develop phase limitations for the discrete-time Zames-Falb multipliers based on the separation theorem for Banach spaces. By contrast with their continuous-time counterparts they lead to numerically efficient results that can be computed either in closed form or via a linear program. The closed-form phase limitations are tight in the sense that we can construct multipliers that meet them with equality. We discuss numerical examples where the limitations are stronger than others in the literature. The numerical results complement searches for multipliers in the literature; they allow us to show, by construction, that the set of plants for which a suitable Zames-Falb multiplier exists is non-convex.
\end{abstract}

\begin{IEEEkeywords}
 Absolute stability; Zames--Falb multipliers; Dual spaces; Integral quadratic constraints. 
\end{IEEEkeywords}

\section{Introduction}

The Zames-Falb multipliers are an important tool for confirming the input-output stability of Lurye systems \cite{Joaquin:2016}. They preserve the positivity of nonlinearities that are either memoryless and monotone or  memoryless, odd  and monotone. Hence, via loop transformation, they may be used for Lurye systems with nonlinearities that are either memoryless and slope restricted or memoryless, odd and slope restricted.

The continuous-time Zames-Falb multipliers were first proposed by O'Shea~\cite{OShea67} and formalized in \cite{Zames:1968}. Similarly their discrete-time counterparts were first proposed by \cite{OShea:1967} and formalized by \cite{Willems:1968,Willems:71}.

There has been considerable interest in Zames--Falb multipliers, in particular, numerical searches in continuous-time domain~\cite{Chen:1996,Turner:2009,Turner:11,Carrasco:2012,Turner2012,Turner:13,Chang:2012,CARRASCO:2014b,turner2019analysis}. Recently, discrete-time domain searches have been presented in~\cite{Joaquin:2019}. Whereas the selection of the parametrisation in continuous-time becomes a fundamental question and remains as an open question~\cite{VEENMAN2016}, FIR multipliers have been shown to be the most effective structure in discrete-time~\cite{Joaquin:2019}. There has also been interest in generalizing the class, both to multi-input, multi-output nonlinearities~\cite{Safonov:00,DAMATO2001,Kulkarni:2002,Mancera:2004} and to nonlinearities outside the original classes considered by Zames and Falb~\cite{Rantzer2001,Materassi2011,Heath2015}. Applications of  the discrete-time Zames--Falb multipliers range from input-constrained model predictive control~\cite{Heath2007,Panos2020} to first order numerical optimization algorithms~\cite{Lessard2016,Michalowsky2020}.

It is useful to establish phase limitations of the available Zames-Falb multipliers. It is trivial to show that their phase must lie between $-90^o$ and $+90^o$. However if this were the only limitation then the Kalman Conjecture would be true for all plants. Instead, fourth-order continuous-time counterexamples \cite{Fitts} and second-order discrete-time counterexamples \cite{Carrasco:2015,Heath:2015} are known. In \cite{Shuai:2018} we develop phase limitations over frequency regions based on an idea by~\cite{Megretski:1995}. It is remarkable that the phase limitations for discrete-time multipliers are more restrictive than those for continuous-time multipliers; in particular in~\cite{Shuai:2018} it is shown that there are phase limitations for discrete-time multipliers over a single frequency interval (with non-zero measure). 

An alternative approach was suggested by J\"onsson and co-workers~\cite{Jonsson:1995,Jonsson:1996,Jonsson:1996a,Jonsson:1996c,Jonsson:1997,Jonsson:1999}.  Specifically J\"onsson develops separation results for continuous-time Zames-Falb multipliers where the nonlinearity is odd. Despite the elegance of the framework, the results have not previously led to useful algorithms (to the authors' knowledge) largely because ``it is in most applications hard to find a suitable frequency grid for the application of the results''. 

In this paper we revisit J\"onsson's approach and develop the discrete-time counterpart for the separation result and its application to Zames--Falb multipliers. The main contribution is to show that the discrete-time results lead to efficient numerical methods and insightful results that are different in kind to the continuous-time results.

%The results significantly improve the results given in \cite{WangTAC} for all the examples considered there, {\color{magenta} [Meaning not quite clear. Do we repeat all examples in \cite{WangTAC}? Do we add some more?]} both in terms of the phase limitations we find and the time taken to compute them.
The results of this paper not only reduce the conservativeness with respect to~\cite{Shuai:2018}, but also significantly reduce the computational burden. In particular for the discrete-time case, phase limitations based on duality results can be computed using a linear program; furthermore there are useful results with a closed-form solution. The closed-form solution phase limitations are tight in the sense that it is possible to construct Zames--Falb multipliers that meet them with equality. 
In the numerical examples we test stability as we vary the loop gain; in all the examples we find only a small gap between gains where we can find Zames-Falb multipliers via the search of~\cite{Joaquin:2019} and gains where the results of this paper show no Zames-Falb multiplier exists.

This has interesting consequences. It allows us to show, by construction, that the set of plants for which a suitable Zames-Falb multiplier exists is non-convex (Theorem~\ref{th:2}). This in turn is strongly indicative of the reason why the Zames-Falb multipliers have not, to-date, been widely used as a design tool. It also adds credence to the conjecture we posed in~\cite{Shuai:2018}; the implications are further explored in~\cite{Seiler:2020}.

The remainder of the paper is structured as follows. In Section II we provide necessary mathematical preliminaries. In Section III we provide the discrete-time counterparts to J\"onsson's results, applied to both odd and nonodd nonlinearities. Section IV represents the main contribution of the paper by developing tractable conditions for duality results. In Section IV.A we discuss results at single frequencies where there is a closed-form solution; in Section IV.B we show these results are tight in the sense that it is possible to construct multipliers that meet the phase limitation with equality; in Section IV.C we consider limitations at several frequencies that can be obtained via linear programming. In Section V we discuss the application to numerical examples, and compare results with both the search for multipliers of \cite{Joaquin:2019} and the phase limitations of \cite{Shuai:2018}.  In two Appendices we provide details about the algorithm used to check conditions in~\cite{Shuai:2018} and provide a complementary continuous-time results to~\cite{Jonsson:1996c} for the case where the nonlinearity is non-odd.

\section{Preliminaries}\label{sec:preliminaries}

\subsection{Mathematical preliminaries}\label{sec:spaces}
Some linear vector spaces, frequency domain spaces and their dual spaces are defined in this part.  

\subsubsection{Linear vector spaces}
Let $\mathds{Z}$ and $\mathds{Z}^+$ ($\mathds{R}$ and $\mathds{R}^+$) be the set of integers (real numbers) and non-negative integers (real numbers), respectively. Let $\mathds{Q}^+$ be the set of non-negative rational numbers.  Let $\mathds{C}$ be the set of complex numbers. 
Two positive integers $a$ and $b$ are said to
be coprime or relatively prime if and only if their greatest
common divisor is $1$, i.e. $gcd(a;b) = 1$, henceforth $a\perp b$.

%Furthermore, the positive cone $\mathbf{P}$ (negative cone $\mathbf{N}$) in a linear vector space $\mathbf{X}$ consists of nonnegative or positive semidefinite  (nonpositive or negative semidefinite) elements in $\mathbf{X}$.    It is clear that    positive (negative) cones are convex.

The space of real (complex) vectors and square matrices is denoted by $\mathds{R}^n$ and $\mathds{R}^{n\times n}$ ($\mathds{C}^n$ and $\mathds{C}^{n\times n}$). For a vector $x\in\mathds{R}^n$, the condition $x\succeq0$ is satisfied if all its elements are non-negative.  

The Hermitian conjugate $A^\sim$ of the complex matrix $A$ is defined as its transpose conjugate. A real matrix is said to be symmetric if $A^\top=A$. A complex matrix is said to the Hermitian if $A^\sim=A$.  The Frobenius norm is defined as $|A|_F=\sqrt{\text{tr}(A^{\sim}A)}$. Let $\mathbf{S}_R^{m\times m}\subset \mathds{R}^{m\times m}$ ($\mathbf{S}_C^{m\times m}\subset \mathds{C}^{m\times m}$) consist of symmetric (Hermitian) matrices with the Frobenius norm. 

A symmetric real matrix $A\in\mathds{R}^{n\times n}$ is said to be positive definite (semidefinite) if $x^\top A x>0$ ($x^\top A x\geq0$) for all non-zero $x\in\mathds{R}^n$. A complex hermitian matrix $A\in\mathds{C}^{n\times n}$ is said to be positive definite (semidefinite) if $\hbox{Re}\{x^\top Ax\}>0$ ($\hbox{Re}\{x^\top Ax\}\geq0$) for all $x\neq0$. 

Let $\mathbf{X}$ and $\mathbf{Y}$ be normed vector spaces. The dual space of $\mathbf{X}$, denoted by $\mathbf{X}^*$, is the Banach space consisting of all bounded linear functionals on $\mathbf{X}$. The real-valued linear functional $\langle x,x^*\rangle$  denotes the value of  $x^*\in \mathbf{X}^*$  at $x\in \mathbf{X}$. %Let $H: \mathbf{X}\mapsto \mathbf{Y}$ be a bounded linear operator. %Then its adjoint operator $H^{\times}: \mathbf{Y}^* \mapsto \mathbf{X}^*$ is defined by the relation $\langle Hx,y^*\rangle=\langle x,H^{\times}y^*\rangle$, for all $x\in \mathbf{X}$ and $y^*\in \mathbf{Y}^*$. Particularly, when $\mathbf{X}$ and $\mathbf{Y}$ are real Hilbert spaces, the operator is said to be self-adjoint if  $H^{\times}=H$. 

A set $\mathbf{C}$ in a linear vector space $\mathbf{X}$ is said to be a cone (with vertex at the origin) if  $x\in \mathbf{C}$ implies that $\alpha x\in \mathbf{C}$,  $\forall \alpha\ge0$. In addition, the cone is convex if  $x_1,x_2\in \mathbf{C}$ implies that $\alpha x_1+\beta x_2\in \mathbf{C}$,  $\forall \alpha,\beta\ge0$.

\subsubsection{Signal spaces}
 Let $\ell$ be the space of all real-valued sequences $h: \mathds{Z}^+\rightarrow\mathds{R}$.  Let $\ell_2$ be the space of all square-summable sequences  $h:\mathds{Z}^+\rightarrow\mathds{R}$. In discrete time, the extended spaces of $\ell_p$ is given by $\ell$, i.e. $\ell_{pe}=\ell$ for all $p$, hence there is no need for the extended notation~\cite{vidyasagar}. Nonetheless, the notation $\ell_{pe}$ is sometimes used in other parts of the literature, e.g.~\cite{Desoer:1975}, to highlight the norm of the original space. Let ${\ell}_1(\mathds{Z})$ be the space of all absolute-summable $h:\mathds{Z}\rightarrow\mathds{R}$. For $h\in\ell_1(\mathds{Z})$, its norm is defined as
\begin{equation}
\|h\|_1=\sum_{i=-\infty}^{\infty}h_i
\end{equation}

The extension of the above definitions to vector-valued functions or sequences is trivial and will be denoted by $\ell^n$, $\ell_2^n$ and $\ell_1^n(\mathds{Z})$.

The discrete-time Fourier transform of $h$ is $\widehat{h}(e^{j\omega})=\sum_{i=-\infty}^{\infty}h_i e^{-j\omega i}$, $\omega\in [-\pi,\pi]$. If $h\in\ell_1$, then $\widehat{h}(e^{j\omega})$ converges for all $\omega$. It can be extended to $\ell_2$ signal by using the limit of the truncated sequence; see~\cite{Korner,oppenheim1999discrete} for further details of the convergence properties. For real sequences, the real part of the Fourier transform is even and the imaginary part is odd. For continuous Fourier transform and its convergence properties, see~\cite{Korner}.

 %Particularly, the dual space of ${\ell}_1$ is ${\ell}_{\infty}$ with the linear functional $\langle x,x^*\rangle=\sum_{i=-\infty}^{\infty}x_i x^*_i$, where $x_i\in {\ell}_1$ and $x^*_i\in{\ell}_{\infty}$.

%Their dual spaces are identified with themselves respectively, and the linear functionals are both real-valued with the  definition $\langle x,x^*\rangle=\text{tr}(xx^*)$.

\subsubsection{System spaces}
Let $\mathbf{RL}_{\infty}^{m\times m}$ be the space consisting of proper real rational transfer function matrices $G:\mathds{C}\rightarrow \mathds{C}^{m\times m}$ that have no pole on the unit circle  in the complex plane.  Let $\mathbf{RH}_{\infty}^{m\times m}$ be a subspace of $\mathbf{RL}_{\infty}^{m\times m}$, where functions have all poles inside the open unit disk. 
%The infinity norm of $G\in \mathbf{RL}_{\infty}^{m\times m}$ is defined as $\|G\|=\sup_{\omega\in [-\pi,\pi]}|G(e^{j\omega})|_F$ \footnote{Note that this is not the standard norm, but an equivalent norm.}.
Let $\mathbf{S}_{\infty}^{m\times m}\subset \mathbf{RL}_{\infty}^{m\times m}$ be the subspace consisting of transfer function matrices satisfying that $G(e^{j\omega})$ is Hermitian, i.e. $G(e^{j\omega})=G(e^{j\omega})^\sim$, $\forall \omega\in [-\pi,\pi]$.

The dual spaces are defined as follows:
\begin{definition}[\cite{Luenberger:1997,Limaye:2016}]\label{df:nbv}
	A function $f:\omega\in[0,\pi] \mapsto \mathbf{S}_C^{m\times m}$ is said to be of bounded variation if its total variation on $[0,\pi]$ is finite, i.e.
	\begin{equation}\label{eq:tv(z)}
	\hbox{T.V.}(f)\triangleq\int_{0}^{\pi}|\mathrm{d} f(\omega)|_F<\infty,
	\end{equation}
	Furthermore, the function is said to be normalised if it
	vanishes at $0$ and is right continuous on $(0,\pi)$. Moreover, $\|f\|=\text{T.V.}(f)$
\end{definition}	

\begin{definition}[Space $S_{NBV}^{m\times m}$~\cite{Jonsson:1996}]\label{df:S_nbv}	
	The space $\mathbf{S}_{NBV}^{m\times m}$ is the Banach space consisting of functions $f:\omega\in[-\pi,\pi] \mapsto \mathbf{S}_C^{m\times m}$ which are normalised with bounded variation on $[0,\pi]$, and which satisfy $f(-\omega)=-\overline{f(\omega)}$  $\forall \omega\in [0,\pi]$.
\end{definition}

Moreover, the convex cone $\mathbf{P}_{NBV}^{m\times m}\subset \mathbf{S}_{NBV}^{m\times m}$ consists of functions that also satisfy $f(\omega_1)\ge f(\omega_2)$, $\forall  \omega_1\ge \omega_2\ge 0$. The underlying methodology uses functions in the real line, hence let us define $g(\omega)=G(e^{j\omega})$ for any given transfer function matrix  $G\in \mathbf{S}_{\infty}^{m\times m}$. 

\begin{lemma}[\cite{Luenberger:1997}]
	The dual of $\mathbf{S}_{\infty}^{m\times m}$ can be identified with $\mathbf{S}_{NBV}^{m\times m}$.  For $G\in \mathbf{S}_{\infty}^{m\times m}$ and $f\in\mathbf{S}_{NBV}^{m\times m}$, the real-valued linear functional is defined by the Stieltjes integral   
	\begin{equation}\label{eq:Stieltjes_integral}
	\langle g ,f\rangle=\int_{-\pi}^{\pi}\text{tr}\left(g(\omega)\mathrm{d} f(\omega)\right).
	\end{equation}
\end{lemma}

\begin{IEEEproof}
	If the function $G\in \mathbf{S}_{\infty}^{m\times m}$  then $g$ belong to the class of continuous functions in the interval $[-\pi,\pi]$. By symmetry, the integral can be rewritten as
	\begin{equation}\label{eq:Stieltjes_integral1}
	\langle g ,f\rangle=2\int_{0}^{\pi}\text{tr}\left(g(\omega)\mathrm{d} f(\omega)\right);
	\end{equation}
	hence the problem can be reduced to the interval $\omega\in [0,\pi]$. Moreover, Definition~\ref{df:S_nbv} in conjunction with the properties of Fourier transform for real sequences ensures that the linear function is real-valued. As a result, the proof is a consequence of the Riesz Representation Theorem (see p. 113 in~\cite{Luenberger:1997}) for the dual of space of continuous functions on a real interval. 
\end{IEEEproof}

\subsubsection{Periodicity of complex exponential functions} 

The periodicity of the exponential function will be exploited in Section~\ref{Sec:application} by using the following result.

\begin{lemma}\label{lemma:periodic}
	Given an $\omega=\frac{\alpha}{\beta}\pi$, where $\alpha\in\mathds{Z}^+$, $\beta \in\mathds{Z}^+$ with $\alpha\perp\beta$ and $\alpha<\beta$, the minimum period $T$ of a complex exponential sequence $h_i=e^{-j\omega i}$, $\forall i\in \mathds{Z}$, is
	\begin{equation}\label{eq:period}
	T=\begin{cases}
	2\beta \quad \text{when }\alpha \text{ is odd},\\
	\beta \quad \text{  when }\alpha \text{ is even}.
	\end{cases}
	\end{equation}
	Moreover, the phase of the complex exponential sequence $h_i=e^{-j\omega i}$ $\forall i\in \mathds{Z}$, is given by the finite sequence 
	\begin{equation}\label{eq:phase}
	\begin{cases}
	\{0,\frac{-\pi}{\beta},\frac{-2\pi}{\beta},...,\frac{-(2\beta-1)\pi}{\beta}\} \quad \text{when }\alpha \text{ is odd},\\
    \{0,\frac{-2\pi}{\beta}, \frac{-4\pi}{\beta},...,\frac{-2(\beta-1)\pi}{\beta}\} \quad \text{  when }\alpha \text{ is even}.
	\end{cases}
	\end{equation}
\end{lemma}

\begin{IEEEproof}
	The identity
	$e^{-j\omega i}=e^{-j\omega(i+T)}=e^{-j\omega i}e^{-j\omega T}$ implies that $e^{-j\omega T}=1$, ie. $\omega T=2 n \pi$, where $n\in \mathds{Z}^+$. Hence, with $\omega=\frac{\alpha}{\beta}\pi$,  $T=\frac{2n\beta}{a}$. Finally, when $\alpha$ is odd, the minimum $T=2\beta$ by setting $n=\alpha$; when $\alpha$ is even, the minimum $T=\beta$ by setting $2n=\alpha$. The phase result then follows trivially.
\end{IEEEproof}

For all other frequencies, i.e. $\omega=\gamma \pi$ with $\gamma$ irrational, we can state the following result:
\begin{lemma}[\cite{Patrick:1995}]~\label{lemma:1}
Let $\omega=\gamma \pi$ with $\gamma\in\mathds{R^+}\backslash \mathds{Q}$. The sequence of complex numbers $e^{-j\omega i}$ for $i=0,1,2,...$ is uniformly dense in the unit circle.	
\end{lemma}

\subsection{Lurye systems}
\begin{figure}[ht]
	\centering
	\ifx\JPicScale\undefined\def\JPicScale{1}\fi
\unitlength \JPicScale mm
\begin{picture}(50,15)(0,0)

\linethickness{0.3mm}
\put(10,2.5){\line(1,0){10}}
\put(10,2.5){\line(0,1){7.5}}
\put(10,10){\vector(0,1){0.12}}
\put(10,12.5){\circle{5}}

%\linethickness{0.5mm}
%\put(12,15){\line(1,0){2}}
\linethickness{0.3mm}
\put(0,12.5){\line(1,0){7.5}}
\put(7.5,12.5){\vector(1,0){0.12}}
\put(12.5,12.5){\line(1,0){7.5}}
\put(20,12.5){\vector(1,0){0.12}}
\put(30,12.5){\line(1,0){10}}
\put(40,5){\line(0,1){7.5}}
\put(40,5){\vector(0,-1){0.12}}

\put(40,2.5){\circle{5}}
\put(42.5,2.5){\line(1,0){7.5}}
\put(42.5,2.5){\vector(-1,0){0.12}}
\put(30,2.5){\line(1,0){7.5}}
\put(30,2.5){\vector(-1,0){0.12}}

\put(20,10){\line(1,0){10}}
\put(20,15){\line(1,0){10}}
\put(20,10){\line(0,1){5}}
\put(30,10){\line(0,1){5}}

\put(20,0){\line(1,0){10}}
\put(20,5){\line(1,0){10}}
\put(20,0){\line(0,1){5}}
\put(30,0){\line(0,1){5}}

\put(25,12.5){\makebox(0,0)[cc]{$G$}}
\put(25,2.5){\makebox(0,0)[cc]{$\phi$}}

\put(3,14){\makebox(0,0)[cc]{$u_1$}}
\put(15,14){\makebox(0,0)[cc]{$e_1$}}
\put(12.5,9){\makebox(0,0)[cc]{$-$}}
\put(47,4){\makebox(0,0)[cc]{$u_2$}}
\put(35,4){\makebox(0,0)[cc]{$e_2$}}
%\put(15,4){\makebox(0,0)[cc]{$h$}}
%\put(35,14){\makebox(0,0)[cc]{$y$}}
\end{picture}
	\caption{Lurye systems}
	\label{fig:lure}
\end{figure}
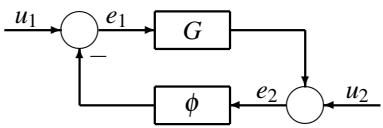

In this paper, we consider SISO Lurye systems in Fig. \ref{fig:lure}, which is expressed as
\begin{equation}\label{eq:Lure}
e_2=u_2+Ge_1, \quad e_1=u_1-\phi(e_2).
\end{equation}
The system (\ref{eq:Lure}) is well-posed if the inverse map $(e_1,e_2) \mapsto (u_1,u_2)$ is causal in $\ell^{2}$. In addition, it is $\ell_2$-stable if it is well-posed, and the signals $(e_1,e_2)$ belong to $\ell_2^2$ for any $(u_1,u_2)\in\ell_2^2$. 

A nonlinear operator  $\phi: \ell \mapsto \ell$ is memoryless if there exists a map $N:  \mathds{R} \to \mathds{R}$ such that $(\phi(\upsilon))_i=N(\upsilon_i)$, $\forall i \in \mathds{Z}$. Assume that $N(0)=0$. 
The memoryless nonlinearity $\phi$ is sector bounded, denoted by $\phi \in [k_1,k_2]$  with $k_1<k_2<\infty$, if $k_1 \le N(x)/x \le k_2, \forall x \in \mathds{R}$. 
The memoryless nonlinearity $\phi$ is monotone if $[N(x_1)-N(x_2)]/(x_1-x_2)\ge 0, \forall x_1, x_2 \in \mathds{R}$ and $x_1 \ne x_2$.
The uncertainty $\phi$ is slope-restricted, denoted by $\phi\in S[k_1,k_2]$, if $k_1(x_1-x_2) \le N(x_1)-N(x_2) \le k_2(x_1-x_2), \forall x_1, x_2 \in \mathds{R}$ and $x_1 \ne x_2$. 
%With a little abuse of the notation, we denote monotone uncertainties as $\phi^{M}\in S[0,\infty]$.
Finally, the  memoryless operator $\phi$ is said to be odd if $N(-x)=-N(x)$, $\forall x\in \mathds{R}$.

\subsection{IQC framework}\label{sec:ZF_multiplier}
The discrete-time extension of the IQC framework is trivial by replacing the imaginary axis with the unit circle.

\begin{definition}[IQC \cite{Megretski:1997,Kao:2012}]\label{df:IQC}	
	Let $\Pi:[0,\pi]\rightarrow \mathbf{S}_{C}^{m\times m}$ be a measurable Hermitian-valued function.  Then, a bounded and causal operator $\phi$ is said to satisfy the IQC defined by $\Pi$, if
	\begin{gather}\label{eq:iqc}
	\int_0^\pi
	\begin{bmatrix}
	\widehat{v}(e^{j\omega}) \\ \widehat{\phi(v)}(e^{j\omega})
	\end{bmatrix}^{\sim}
	\Pi(e^{j\omega})
	\begin{bmatrix}
	\widehat{v}(e^{j\omega}) \\ \widehat{\phi(v)}(e^{j\omega})
	\end{bmatrix}
	\mathrm{d}\omega \ge 0, \quad \forall v\in\ell_2,
	\end{gather}
	where $\widehat{v}$ and  $\widehat{\phi(v)}$ denote the  discrete Fourier transform of $v$ and $\phi(v)$ respectively.
\end{definition}

\begin{remark}
	It is standard in the IQC literature to restrict the definition of the IQC to the interval $[0,\pi]$, but it can be extended to the interval $[-\pi,\pi]$, which is the natural interval for duality results, when $\Pi$ is constructed as the z-transform of a time-domain operator.
\end{remark}

Similarly, we can write the IQC theorem where we restrict our attention to continuous multipliers: 
\begin{theorem} [\cite{Megretski:1997}] \label{th:IQC_discrete}
	Let $\Pi:[0,\pi]\rightarrow \mathbf{S}_{C}^{m\times m}$ a continuous measurable Hermitian-valued function. For the system in Fig. \ref{fig:lure}, let $G \in \mathbf{RH}_{\infty}$  and $\phi$ be a causal bounded operator. Assume that $\forall \tau \in [0,1]$,
	\begin{enumerate}
		\item  the feedback interconnection between $G$ and $\tau \phi$ is well-posed;
		\item the operator $\tau \phi$ satisfies the IQC defined by $\Pi$;
		\item the following frequency domain inequality holds,
		\begin{gather}\label{eq:iqc_stability}
		\begin{bmatrix}
		G(e^{j\omega})  \\
		1 \\
		\end{bmatrix}^{\sim}
		\Pi(e^{j\omega})
		\begin{bmatrix}
		G(e^{j\omega})  \\
		1 \\
		\end{bmatrix}
		> 0, \;\qquad \forall \omega\in[0,\pi].
		\end{gather}
	\end{enumerate}
	Then, the system in Fig.\ref{fig:lure} is $\ell_2$-stable.
\end{theorem}

\begin{remark}
	The corresponding condition in the classical theory (e.g. \eqref{eq:iqc_stability} in Megretski and Rantzer~\cite{Megretski:1997}) has an additional term $\epsilon I$ on the right hand side. In discrete-time, the continuity of the multiplier, which is required to define the dual spaces, allow us to remove this condition as the condition is tested in a bounded interval.
\end{remark}

\subsection{Zames-Falb multipliers}

If $\phi$ is memoryless, bounded and monotone, then the Zames-Falb IQC with the class of Zames-Falb multipliers is valid to preserve the positivity in (\ref{eq:iqc}). A wider class of Zames-Falb multipliers is obtained when $\phi$ is also odd.

\begin{definition}[Discrete-time LTI Zames-Falb multipliers~\cite{OShea:1967,Willems:1968}]
	Let the real-valued sequence $h$ and the operator $H$ be a discrete-time Fourier transform pair. Then, an operator $M(e^{j\omega})=1-H(e^{j\omega})$  belongs to the set $\mathcal{M}$ if $h_0=0$, $h_i \ge 0$ and  $ \|h\|_1\leq1$. Similarly, the operator $M$ belongs to the set $\mathcal{M}_{\text{odd}}$ if
	$h_0=0$ and $\|h\|_1\leq1$.
\end{definition}

\begin{remark}
	Although we may use a strictly inequality from a practical point of view,  we use the original definition with non-strict inequalities as the subset of multipliers with $\|h\|_1=1$ will be important for some parts of the paper.
\end{remark}

\begin{remark}
Note that $\mathcal{M}\subset\mathcal{M_{\text{odd}}}$.
\end{remark}

 We use the definition which is analogous to the continuous Zames-Falb multipliers, see~\cite{Zames:1968,Joaquin:2016}. However, the classical definition of discrete Zames-Falb multipliers~\cite{Willems:1968} also includes LTV multipliers, though they have never been used to the best of the authors'  knowledge.

\begin{definition}[Zames-Falb IQC multipliers]\label{df:convex_cone}
	The convex cone $\Pi_{\phi}\subset\mathbf{S}_{\infty}^{2\times 2}$ consists of function matrices $\Pi$, which is in the form
	\begin{equation}\label{eq:zf_iqc}
	\Pi(e^{j\omega})=
	\begin{bmatrix} 0 &  M(e^{j\omega})^{\sim} \\ M(e^{j\omega}) & 0  \end{bmatrix},
	\end{equation}
	with $M\in\mathcal{M}$ or $M\in\mathcal{M_{\text{odd}}}$.
\end{definition}

Then we have the following corollary of Theorem~\ref{th:IQC_discrete} for bounded and monotone nonlinearities.

\begin{corollary}\label{corollary:iqc_stability1}
For the system in Fig. \ref{fig:lure}, let $G \in \mathbf{RH}_{\infty}$,  and let $\phi$ be a bounded and monotone memoryless nonlinearity. The system is $\ell_2$-stable, if there exists  $\Pi\in\Pi_{\phi}$ with $M\in\mathcal{M}$, such that 
\begin{gather}\label{eq:iqc_stability1}
\begin{bmatrix}
G(e^{j\omega})  \\
1 \\
\end{bmatrix}^{\sim}
\Pi(e^{j\omega})
\begin{bmatrix}
G(e^{j\omega})   \\
1 \\
\end{bmatrix}
> 0,\quad  \forall\omega\in[0,\pi].
\end{gather}
Moreover, $M\in\mathcal{M_{\text{odd}}}$ in $\Pi$ if $\phi$ is also odd.
\end{corollary}

\subsection{Phase limitation of Zames-Falb multipliers in~\cite{Shuai:2018}}\label{sec:phase_limit}
In this part, we repeat the phase properties of Zames-Falb multipliers in~\cite{Shuai:2018} in order to compare with the duality approach in this paper. 
\begin{definition}
	Let $0\le a<b\le \pi$. 	Define
	\begin{equation*}
	\psi^d(n)=\frac{\cos(an)-\cos(bn)}{n},\; \phi^d(n)=\frac{\sin(an)-\sin(bn)}{n}.
	\end{equation*}
	
	Then, 
	\begin{subequations}\label{eq:phase_limit}
		\begin{align}
		&\mu^d(n)=\frac{|\psi^d(n)|}{(b-a)+\phi^d(n)},\quad  \rho^d=\max_{n\in \mathds{Z}^+}\mu^d,\\
		&\mu^d_{odd}(n)=\frac{|\psi^d(n)|}{(b-a)-|\phi^d(n)|},\quad \rho^d_{odd}=\max_{n\in \mathds{Z}^+}\mu^d_{odd},
		\end{align}
	\end{subequations}
	where both $\rho^d$ and $\rho^d_{odd}$ are positive and well-defined ($\rho^d,\rho^d_{odd}<\infty$). 
\end{definition}

\begin{theorem}[Phase limitation of Zames-Falb multipliers~\cite{Shuai:2018}]\label{th:phase_limitation} For a discrete-time Zames-Falb multiplier $M$, if there exist $\rho>0$ such that
	\begin{equation}
	Im\left\{M(e^{j\omega})\right\} > \rho \hbox{Re}\left\{M(e^{j\omega})\right\},\quad \forall \omega\in[a,b],
	\end{equation}
	then $\rho<\rho^d$ if $M\in\mathcal{M}$; $\rho<\rho^d_{odd}$ if $M\in\mathcal{M}_{\text{odd}}$.
	
\end{theorem}

\begin{remark}
	With a given frequency pair $a$ and $b$, the operators $\mu^d$ and $\mu^d_{odd}$ $\to 0$ with  $n\to \infty$. Hence, $\rho^d$ and $\rho^d_{odd}$ can be obtained by searching in the range $n\in[1,N]$ with some $N$ being sufficiently large.
\end{remark}

\section{Discrete-time duality results}
%Refer to \cite{Zames:1968} for the definition of continuous-time Zames-Falb multiplier and \cite{Joaquin:2016} for a tutorial. 

\subsection{General separation result and its appplication to Zames--Falb multipliers}

This paper builds on the separation result developed by J\"onsson in the continuous domain (Theorem~4.2 in~\cite{Jonsson:1996}). We restrict our attention to its discrete-time counterpart for SISO plants although extensions to square MIMO plants are straightforward.

\begin{theorem}\label{th:1}
	Let $\Pi_{\phi}\subset\mathbf{S}_{\infty}^{2\times 2}$ be a convex set and let $G\in\mathbf{RH}_\infty$. The following two statements are equivalent:
	\begin{itemize}
		\item There is no $\Pi\in\Pi_\phi$ such that 
		\begin{equation}
		\begin{bmatrix}
			G(e^{j\omega})  \\
			1 \\
		\end{bmatrix}^{\sim}
		\Pi(e^{j\omega})
		\begin{bmatrix}
			G(e^{j\omega})   \\
			1 \\
		\end{bmatrix}
		> 0,\qquad  \forall \omega\in[0,\pi].
		\end{equation}
		\item There exits a nonzero $Z\in\mathbf{P}_{NBV}$ such that
	\begin{equation}
\int_{-\pi}^\pi
\text{tr}\left(\Pi(e^{j\omega})\begin{bmatrix}
G(e^{j\omega})  \\
1 \\
\end{bmatrix}
\mathrm{d}Z(\omega)
\begin{bmatrix}
G(e^{j\omega})   \\
1 \\
\end{bmatrix}^{\sim}\right)
\leq 0,
\end{equation}
for all $\Pi\in\Pi_\phi$.		
\end{itemize}
\end{theorem}
 
\begin{IEEEproof}Although the proof given in~\cite{Jonsson:1996} is for continuous time only, it is based on the separation principle for Banach spaces and the properties of linear operators. It can thus be directly translated to discrete time.
\end{IEEEproof}

As we will restrict our attention SISO Zames-Falb multipliers, Theorem~\ref{th:1} can be rewriten as follows:

\begin{corollary}\label{cor:10}
		Let $G\in\mathbf{RH}_\infty$. The following two statements are equivalent:
	\begin{itemize}
		\item There is no $M\in\mathcal{M}$ (or $M\in\mathcal{M}_\text{odd}$) such that 
		\begin{equation}
		\hbox{Re}(M(e^{j\omega})G(e^{j\omega}))> 0,\qquad  \forall \omega\in[0,\pi].
		\end{equation}
		\item There exits a nonzero $Z\in\mathbf{P}_{NBV}$ such that
		\begin{equation}
		\int_{-\pi}^\pi\hbox{Re}(M(e^{j\omega})G(e^{j\omega}))\mathrm{d}Z(\omega)
		\leq 0.
		\end{equation}
		for all $M\in\mathcal{M}$ (or $M\in\mathcal{M}_\text{odd}$).		
	\end{itemize}		
\end{corollary}

\begin{IEEEproof}
It follows by applying Theorem~\ref{th:1} and substituting $\Pi_\phi$ with the set Zames-Falb IQC multiplier in Definition~\ref{df:convex_cone}. Note that the set is convex by definition. 
\end{IEEEproof}	

From a computational point of view, it is interesting to use a finite parametrization of $Z$ by using the atomic measure space, i.e. for $0<\omega_{1}<\omega_{2}<\dots<\omega_{N}\leq\pi$ and non-negative $\lambda_1,\lambda_2,\dots,\lambda_N$, let us define
\begin{equation}\label{eq:zw_discrete}
Z(\omega)=\sum_{r=1}^{N}\left\{ \lambda_r\theta(\omega-\omega_r) -  \lambda_r\theta(-\omega-\omega_r)\right\},\quad \omega\in[-\pi,\pi],
\end{equation}
where $\theta$ is the step function. It is trivial that $Z\in\mathbf{P}_\text{NBV}$.

\begin{corollary}\label{cor:1} Let $G\in\mathbf{RH}_\infty$. Assume there exist $0<\omega_{1}<\omega_{2}<\dots<\omega_{N}\leq\pi$ and non-negative $\lambda_1,\lambda_2,\dots,\lambda_N$, where at least one $\lambda_r$ is nonzero, such that
	\begin{equation}\label{eq:2}
 \sum_{r=1}^N \hbox{Re}\{\lambda_r M(e^{j\omega_{r}}) G(e^{j\omega_{r}})\}\leq 0,
	\end{equation}		
	for all $M\in\mathcal{M}$ ($M\in\mathcal{M_{\text{odd}}}$), then there is no $M\in\mathcal{M}$ ($M\in\mathcal{M_{\text{odd}}}$) such that 
	\begin{equation}
	\hbox{Re}\{M(e^{j\omega}) G(e^{j\omega})\}> 0,\quad  \forall \omega\in[0,\pi].
	\end{equation}	
\end{corollary}
\begin{IEEEproof}
	The proof follow from Corollary~\ref{cor:10} by using the parametrization of $Z(\omega)$ as the atomic measure in~\eqref{eq:zw_discrete}.
	\end{IEEEproof}

Loosely speaking, we will be able to ensure that there is no suitable Zames-Falb multiplier for a given plant $G$, if we are able to show that all Zames-Falb multipliers satisfy \eqref{eq:2}. The rest of the paper deals with the development of conditions on $G$ such that we can ensure that~\eqref{eq:2} is satisfied for all Zames-Falb multipliers.

\subsection{Duality conditions}\label{sec:slope_restricted_nonlinearity}

We now develop conditions over $G$ ensuring that the duality result in Corollary~\ref{cor:1} is satisfied for all discrete LTI Zames-Falb multipliers.

Before providing the main result, the following lemma shall be used during the proof of the main result:
\begin{lemma} For any given $0 < \omega_1\le \cdots \le \omega_N \le \pi$, and $\lambda_1, \cdots, \lambda_N \ge 0$. Then
	\begin{equation}
		\min_{i\in \mathds{Z}}\left[\sum_{r=1}^{N} \hbox{Re}\left\{\lambda_1 G(e^{j\omega_r})e^{-j\omega_r i}\right\}\right]\leq 0.
	\end{equation}
\end{lemma}
\begin{IEEEproof} Initially, we prove for the case $N=1$. If $\lambda_1=0$, the results is trivially true. If $\lambda_1>0$,
	Let us assume that 
	\begin{equation}\label{eq:12}
\xi_i^{\omega_1}=\hbox{Re}\left\{\lambda_1 G(e^{j\omega_1})e^{-j\omega_1 i}\right\}>0
\end{equation}
for all $i\in\mathds{Z}$, then $\sum_{i=0}^{\infty} \xi_i^{\omega_{1}}>0$

The sequence $\xi_i^{\omega_{1}}$ can be either periodic or aperiodic.

If $\xi_i^{\omega_{1}}$ is periodic with period $T$, then
	\begin{multline}\label{eq:14}
\sum_{i=0}^{T-1}\xi_i^{\omega_{1}}=\sum_{i=0}^{T-1} \hbox{Re}\left\{\lambda_1 G(e^{j\omega_1})e^{-j\omega_1 i}\right\}=\\ \hbox{Re}\left\{\lambda_1 G(e^{j\omega_1})\sum_{i=0}^{T-1}e^{-j\omega_1 i}\right\}
\end{multline}
From the periodicity, any $\omega_1=n_1\frac{2\pi}{T}$ for $n_1\in\mathds{N^+}$, hence it holds
\begin{equation}\label{eq:15}
\sum_{i=0}^{T-1}e^{-jn_1\frac{2\pi}{T}i}=0
\end{equation}
as it can be seen as the second component of the discrete-time Fourier transform of a vector of ones. As a result, $\sum_{i=0}^{T-1} \xi_i=0$, hence the result is obtained by contradiction.

If $\xi_i^{\omega_{1}}$ is aperiodic, for any $\omega_1$,
\begin{equation}\label{eq:16}
\sum_{i=0}^{\infty}e^{-j\omega_1 i}=0
\end{equation}
since by Lemma~\ref{lemma:1} the sequence $e^{-j\omega_1 i}$ is uniformly dense in the unit circle. As a result, $\sum_{i=0}^{\infty} \xi_i^{\omega_{1}}=0$, hence the result is obtained by contradiction.

For $N>1$, the result follows by applying the above argument to every $\xi^{\omega_{r}}_i$ for $r=1,\dots,N$. 

\end{IEEEproof}
\begin{theorem} \label{th:monotone}
Let $G\in\mathbf{RH}_\infty$. Assume there exist $0 < \omega_1\le \cdots \le \omega_N \le \pi$, and $\lambda_1, \cdots, \lambda_N \ge 0$, where at least one $\lambda_r$ is nonzero. If
		\begin{equation}\label{eq:1}
		\sum_{r=1}^{N} \hbox{Re}\left\{\lambda_rG(e^{j\omega_r})\right\}
		\le 
		\min_{i\in \mathds{Z}}\left[\sum_{r=1}^{N} \hbox{Re}\left\{\lambda_r G(e^{j\omega_r})e^{-j\omega_r i}\right\}\right],
		\end{equation}
	then there is no Zames-Falb multiplier $M\in\mathcal{M}$ such that
		\begin{equation}
		\hbox{Re}\left\{M(e^{j\omega})G(e^{j\omega})\right\}>0, \quad \forall \omega\in[0,\pi]. 
		\end{equation}	
	Similarly, if 
		\begin{multline}\label{eq:13}
		\sum_{r=1}^{N} \hbox{Re}\left\{\lambda_r G(e^{j\omega_r})\right\}
		\le 
		-\max_{i\in \mathds{Z}}\left|\sum_{r=1}^{N} \hbox{Re}\left\{\lambda_r G(e^{j\omega_r})e^{-j\omega_r i}\right\}\right|,
		\end{multline}
	then there is no Zames-Falb multiplier $M\in\mathcal{M}_{odd}$ such that
	\begin{equation}
	\hbox{Re}\left\{M(e^{j\omega})G(e^{j\omega})\right\}>0, \quad \forall \omega\in[0,\pi]. 
	\end{equation} 
\end{theorem}
\begin{IEEEproof} The proof follows from Corollary~\ref{cor:1}. We need to show that~\eqref{eq:1} implies~\eqref{eq:2} for all $M\in\mathcal{M}$ or $M\in\mathcal{M_{\text{odd}}}$.

Firstly, let us consider $M\in\mathcal{M}$. We can include the definition of Zames--Falb multipliers as follows. The condition
\begin{equation}\label{eq:4}
\sum_{r=1}^{N} \hbox{Re}\left\{\lambda_rG(e^{j\omega_r})\right\}
\le \sum_{r=1}^{N} \hbox{Re}\left\{\lambda_r G(e^{j\omega_r})e^{-j\omega_r i}\right\},
\end{equation}
for all $i\in\mathds{Z}$ implies that the left hand side must be non-positive. Then, for all sequences $h$ such that $h_i\geq0$, $h_0=0$, and $\|h\|_1\leq1$, it follows that
\begin{equation}\label{eq:5}
\left(1-\sum_{i=-\infty}^\infty h_i\right)\sum_{r=1}^{N} \hbox{Re}\left\{\lambda_rG(e^{j\omega_r})\right\}\leq0.
\end{equation}
Hence substituting~\eqref{eq:4} in~\eqref{eq:5} for the each $i$ in the summation, yields
\begin{multline}
\sum_{r=1}^{N} \hbox{Re}\left\{\lambda_rG(e^{j\omega_r})\right\}
-\\ \sum_{i=-\infty}^{\infty}h_i\left[\sum_{r=1}^{N} \hbox{Re}\left\{\lambda_r G(e^{j\omega_r})e^{-j\omega_r i}\right\}\right]\leq0,
\end{multline}
for all sequences $h$ such that $h_i\geq$, $h_0=0$, and $\|h\|_1\leq1$. Rearranging gives
\begin{equation}
\sum_{r=1}^{N} \hbox{Re}\left\{\lambda_r \left(1-\sum_{i=-\infty}^{\infty}h_i e^{-j\omega_r i}\right)G(e^{j\omega_r})\right\}\leq0.
\end{equation} 
Hence
\begin{equation}
\sum_{r=1}^{N} \hbox{Re}\left\{\lambda_r M(e^{j\omega_{r}})G(e^{j\omega_r})\right\}\leq0,
\end{equation} 	
for all $M\in\mathcal{M}$. Then the result is obtained by using Corollary~\ref{cor:1}.

Following the same approach, let us now consider $M\in\mathcal{M_{\text{odd}}}$. In this case, we have 
\begin{equation}\label{eq:7}
\sum_{r=1}^{N} \hbox{Re}\left\{\lambda_r G(e^{j\omega_r})\right\}
\le 
-\sum_{r=1}^{N} \hbox{Re}\left\{\lambda_r G(e^{j\omega_r})e^{-j\omega_r i}\right\},
\end{equation}
and
\begin{equation}\label{eq:8}
\sum_{r=1}^{N} \hbox{Re}\left\{\lambda_r G(e^{j\omega_r})\right\}
\le 
\sum_{r=1}^{N} \hbox{Re}\left\{\lambda_r G(e^{j\omega_r})e^{-j\omega_r i}\right\},
\end{equation}
for all $i$, and as previously it is trivial that the left hand side must be non-positive.
 Then, for all sequences $h$ such that $h_0=0$ and $\|h\|_1\leq1$, it follows that 
\begin{equation}\label{eq:6}
\left(1-\sum_{i=-\infty}^\infty |h_i|\right)\sum_{r=1}^{N} \hbox{Re}\left\{\lambda_rG(e^{j\omega_r})\right\}\leq0.
\end{equation}
Hence substituting~\eqref{eq:7} (if $h_i<0$) or~\eqref{eq:8} (if $h_i>0$) in~\eqref{eq:6} for each $i$ in the summation, yields
\begin{multline}
\sum_{r=1}^{N} \hbox{Re}\left\{\lambda_rG(e^{j\omega_r})\right\}
-\\ \sum_{i=-\infty}^{\infty}h_i\left[\sum_{r=1}^{N} \hbox{Re}\left\{\lambda_r G(e^{j\omega_r})e^{-j\omega_r i}\right\}\right]\leq0,
\end{multline} 
for all sequences $h$ such that $h_0=0$ and $\|h\|_1\leq1$. Rearranging gives
\begin{equation}
\sum_{r=1}^{N} \hbox{Re}\left\{\lambda_r \left(1-\sum_{i=-\infty}^{\infty}h_i e^{-j\omega_r i}\right)G(e^{j\omega_r})\right\}\leq0,
\end{equation} 
Hence
\begin{equation}
\sum_{r=1}^{N} \hbox{Re}\left\{\lambda_r M(e^{j\omega_{r}})G(e^{j\omega_r})\right\}\leq0,
\end{equation} 	
for all $M\in\mathcal{M}_\text{odd}$. Then the result is obtained by using Corollary~\ref{cor:1}.
\end{IEEEproof}

\section{Main results}
\label{Sec:application}
 The results developed in the previous section can be considered the discrete-time counterpart of results in~\cite{Jonsson:1996} (and Appendix B). The main contribution of this paper is the development of computationally tractable conditions the duality results of Zames-Falb multipliers, which do not exist in continuous-time. 
 
 In particular, we will show that a closed-form expression for limitations at a single frequency is possible. These results represent a change in the understanding of the stability of SISO Lurye system; as they can be easily used by any user by using a Bode Plot. It is possible to improve the single frequency result by including more frequencies, where a linear program is required to derive the phase limitations.

%In this section, firstly we apply Corollary \ref{th:sloperes} for the case $N=1$ as then a closed-form solution is obtained and a phase limitation result is obtained. Secondly, we provide a linear program to efficiently solve the case $N>1$. In some examples, the closed-form solution for $N=1$ can improve the linear program result as it does not require any numerical algorithm. In addition, we show that the case $N=1$ improves the phase limitation results in~\cite{Shuai:2018}, and it can be translated into a phase limitation condition, generating a powerful graphical interpretation of the results. 

\subsection{Single frequency condition as a phase limitation of multipliers}
First, we rewrite the duality conditions in the single frequency case directly.

\begin{corollary}\label{cor:single}
	Let $G\in\mathbf{RH}_{\infty}$. If there exists $\omega_1\in(0,\pi]$ such that 
		\begin{equation}\label{eq:one_f_w}
		 \hbox{Re}\left\{G(e^{j\omega_1})(1-e^{-j\omega_1 i})\right\}\le 0
		, \quad \forall i\in\mathds{Z}
		\end{equation}
	then there is no $M\in\mathcal{M}$ such that
	\begin{equation}
	\hbox{Re}\left\{M(e^{j\omega})G(e^{j\omega})\right\}>0, \quad \forall \omega\in[0,\pi]. 
	\end{equation}
	Similarly, if there exists $\omega_1\in(0,\pi]$ such that 
			\begin{equation}\label{eq:one_f_odd_w}
 \hbox{Re}\left\{G(e^{j\omega_1}) (1\pm e^{-j\omega_1 i})\right\}\le 0, \quad \forall i\in\mathds{Z}
			\end{equation} 
then there is no $M\in\mathcal{M}_{\text{odd}}$ such that
\begin{equation}
\hbox{Re}\left\{M(e^{j\omega})G(e^{j\omega})\right\}>0, \quad \forall \omega\in[0,\pi]. 
\end{equation}		
\end{corollary}
\begin{IEEEproof}
	The result follows from Theorem~\ref{th:monotone} for the case $N=1$ by taking $\lambda_1=1$ without loss of generality. The nonodd case is straightforward, while in the odd case, the absolute value in (\ref{eq:13}) implies that both positive and negative signs of $e^{-j\omega_1 i}$ are in (\ref{eq:one_f_odd_w}).
\end{IEEEproof}	

\begin{remark}
	If $G(e^{j\omega_1})$ is real, the result leads to a trivial conclusion just depending on the sign of $G(e^{j\omega_1})$, i.e. conditions \eqref{eq:one_f_w} amd \eqref{eq:one_f_odd_w} are true if and only if $G(e^{j\omega_1})<0$. If $G(e^{j\omega_1})$ is complex, these conditions can only hold if $\hbox{Re}\left\{G(e^{j\omega_1})\right\}\le 0$, which can be seen a natural consequence of the Circle Criterion.
\end{remark}

\begin{remark}
	Although similar result can be found in continuous-time, it is straightforward to show that the result will only hold for the case $G(s)\in\mathds{R}$ and $G(s)\leq0$. As a result, it will not provide any novel information as this condition can be taken without loss of generality, see~\cite{Joaquin:2013,Carrasco:14a} for further details.
\end{remark}

We have used slightly different notation for the conditions in Corollary~\ref{cor:single} as it provides the following interpretation: At any given frequency, if none of the FIR multiplier with single term at the boundary, i.e. $1-e^{j\omega i}$ for all $i$, is able to correct the lack of positivity of $G$, then no other Zames-Falb multiplier in $\mathcal{M}$ will be able to correct it. For $\mathcal{M_{\text{odd}}}$, we need to test all $1\pm e^{j\omega i}$ to ensure that no other Zames-Falb multiplier in $\mathcal{M_{\text{odd}}}$ will be able to correct. 
	
The result can be rewritten in terms of a phase limitation as in~\cite{Shuai:2018}.
\begin{theorem}\label{Th:PL}
	Let $\alpha\in\mathds{Z^+},\beta\in\mathds{Z^+}$ with $\alpha\perp \beta$ and $\alpha<\beta$. Take $\omega_{1}=\frac{\alpha}{\beta}\pi$. Then 
\begin{eqnarray}
\label{eq:9a}\left|\angle M(e^{j\omega_1})\right|\leq\frac{\pi}{2}\left(1-\frac{1}{\beta}\right),\quad \text{ when $\alpha$ is odd,}\\
\label{eq:9b}\left|\angle M(e^{j\omega_1})\right|\leq \frac{\pi}{2}\left(1-\frac{2}{\beta}\right),\quad \text{ when $\alpha$ is even,}
\end{eqnarray}
for all $M\in\mathcal{M}$ with $M(e^{j\frac{\alpha}{\beta}\pi})\neq 0$. 
Similarly,
\begin{equation}
\left|\angle M(e^{j\frac{\alpha}{\beta}\pi})\right|\leq\frac{\pi}{2}\left(1-\frac{1}{\beta}\right),
\end{equation}
for all $M\in\mathcal{M}_{\text{odd}}$ with $M(e^{j\omega_{1}})\neq 0$.
Moreover, $\angle M(e^{j\pi})=0$ if $M(e^{j\pi})\neq 0$ for all $M\in\mathcal{M}_{\text{odd}}$.
\end{theorem}

\begin{figure*}[p]
	\centering
	\includegraphics[width=\linewidth]{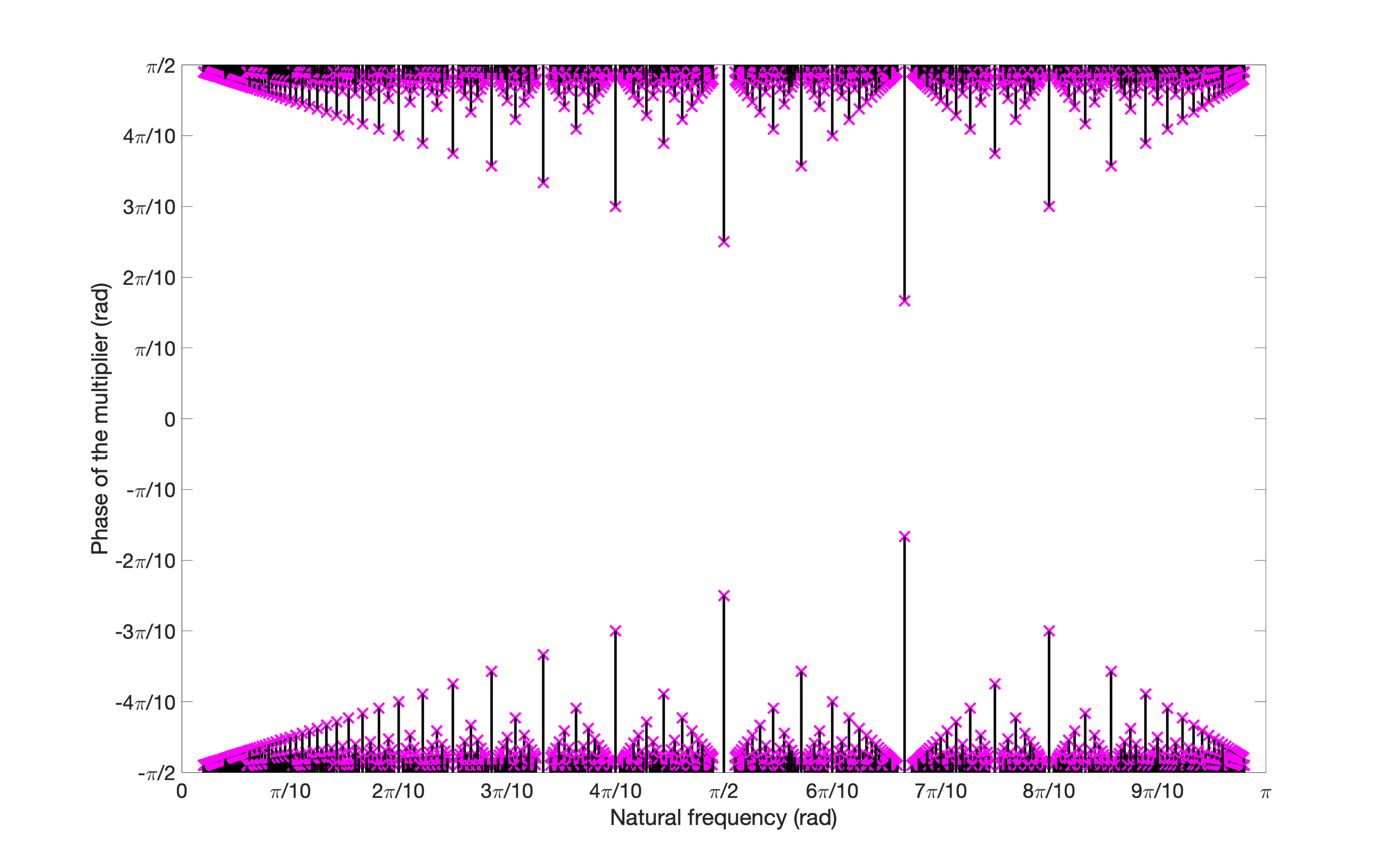}
	\caption{Phase limitation of the class $\mathcal{M}$ for frequencies $\omega=\frac{\alpha}{\beta}\pi$ with $\beta\leq 50$.}
	\label{fig:pl}
	\includegraphics[width=\linewidth]{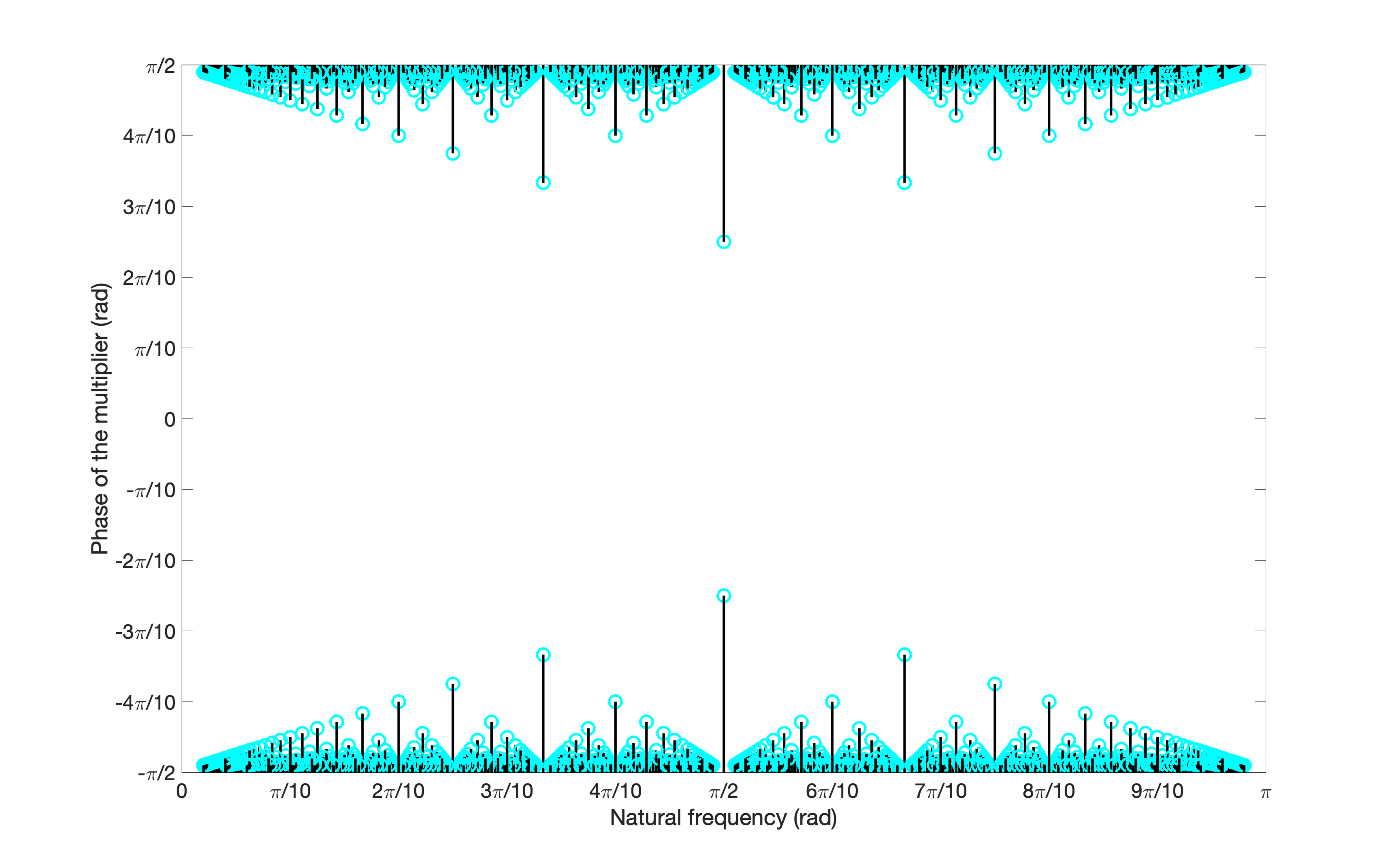}
	\caption{Phase limitation of the class $\mathcal{M}_\text{odd}$ for frequencies $\omega=\frac{\alpha}{\beta}\pi$ with $\beta\leq 50$.}
	\label{fig:pl_odd}
\end{figure*}

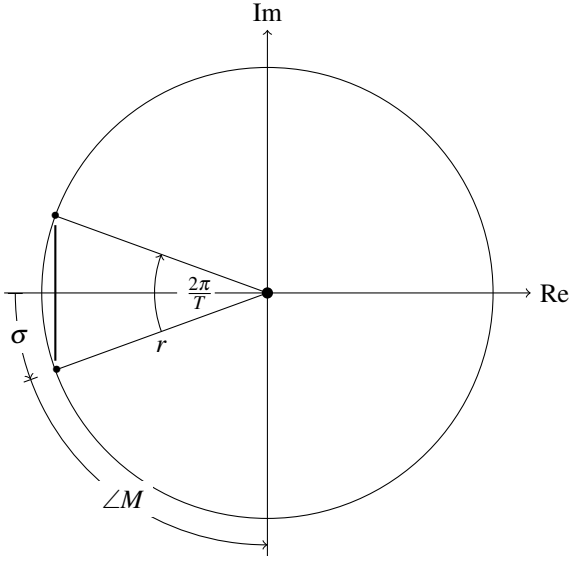
\begin{figure}
	\begin{tikzpicture}
		\draw[-] (-3.5,0) -- (-1.2,0) node[right] {};
	\draw[->] (-0.7,0) -- (3.5,0) node[right] {$\hbox{Re}$};
	\draw[->] (0,-3.5) -- (0,3.5) node[above] {$\hbox{Im}$};
	% the origin
	\coordinate (O) at (0,0);
	% the circle and the dot at the origin
	\draw (O) node[circle,inner sep=1.5pt,fill] {} circle [radius=3cm];
	% the ``\theta'' arc
\foreach \Point in {(160.5:2.99cm), (200.5:2.99cm)}{
	\node at \Point {\textbullet};
}
	\draw 
	(200:3cm) coordinate (xcoord) -- 
	node[midway,below] {$r$} (O) -- 
	(160:3cm) coordinate (slcoord)
	pic [draw,<-,angle radius=1.5cm,"$\frac{2\pi}{T}$"] {angle = slcoord--O--xcoord};
	\node (a1) at (200:3cm) {};
	\node (a2) at (160:3cm) {}; 
    \draw[thick] (a1) to node {} (a2);
	% the outer ``s'' arc
	\draw[|->]
	(-3cm-10pt,0)
	arc[start angle=180,end angle=200,radius=3cm+10pt]
	node[midway,fill=white] {$\sigma$};
	\draw[|->]
	(200:3cm+10pt)
	arc[start angle=200,end angle=270,radius=3cm+10pt]
	node[midway,fill=white] {$\angle M$};
	\end{tikzpicture}
	\caption{The maximum allowed phase for $G(e^{j\omega_1})$, i.e. $-\pi+\sigma$, to ensure that there is no suitable multiplier $M\in\mathcal{M}$ implies in turn a maximum phase for the class of multiplier to recover the positivity of $\hbox{Re}\{M(e^{j\omega_1})G(e^{j\omega_1})\}$. The limitation depends on the period $T$.}\label{fig:1}
\end{figure} 
\begin{figure}
	\begin{tikzpicture}
		\draw[-] (-3.5,0) -- (-1.2,0) node[right] {};
\draw[->] (-0.7,0) -- (3.5,0) node[right] {$\hbox{Re}$};
	\draw[->] (0,-3.5) -- (0,3.5) node[above] {$\hbox{Im}$};
	% the origin
	\coordinate (O) at (0,0);
	% the circle and the dot at the origin
	\draw (O) node[circle,inner sep=1.5pt,fill] {} circle [radius=3cm];
	% the ``\theta'' arc
	\foreach \Point in {(160.5:2.99cm), (200.5:2.99cm)}{
		\node at \Point {\textbullet};
	}
	\draw 
	(200:3cm) coordinate (xcoord) -- 
	node[midway,below] {$r$} (O) -- 
	(160:3cm) coordinate (slcoord)
	pic [draw,<-,angle radius=1.5cm,"$\frac{2\pi}{T}$"] {angle = slcoord--O--xcoord};
	\node (a1) at (200:3cm) {};
	\node (a2) at (160:3cm) {}; 
	\draw[thick] (a1) to node {} (a2);
	% the outer ``s'' arc
	\draw[|->]
	(-3cm-10pt,0)
	arc[start angle=180,end angle=160,radius=3cm+10pt]
	node[midway,fill=white] {$\sigma$};
	\draw[|->]
	(160:3cm+10pt)
	arc[start angle=160,end angle=90,radius=3cm+10pt]
	node[midway,fill=white] {$\angle M$};
	\end{tikzpicture}
	\caption{The minimum allowed phase for $G(e^{j\omega_1})$, i.e. $\pi-\sigma$, to ensure that there is no suitable multiplier $M\in\mathcal{M}$ implies in turn a minimum phase for the class of multiplier to recover the positivity of $\hbox{Re}\{M(e^{j\omega_1})G(e^{j\omega_1})\}$.}\label{fig:3}
\end{figure}
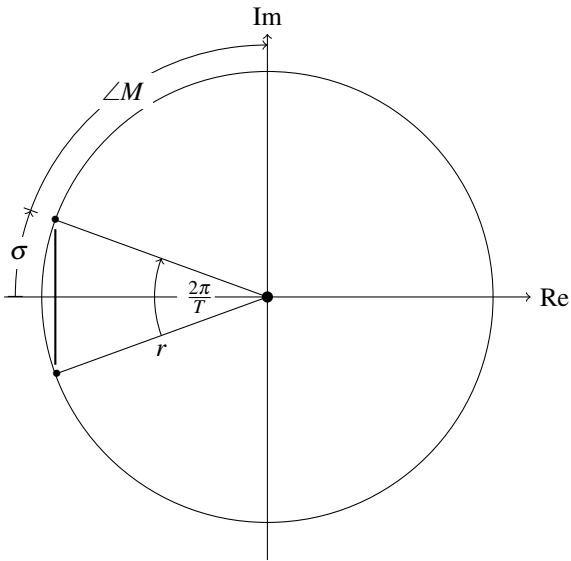

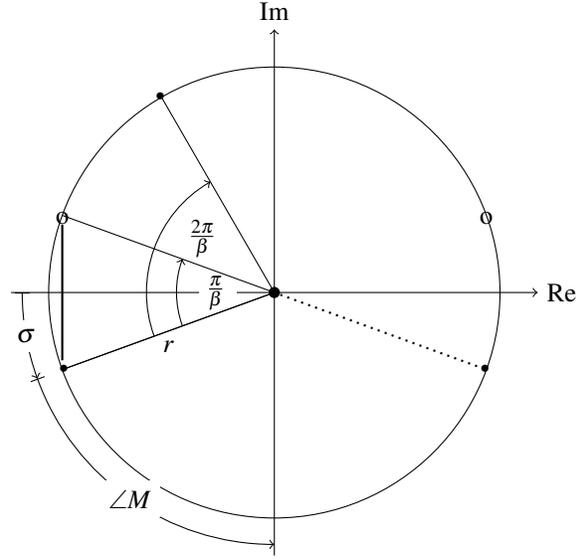
\begin{figure}
	\begin{tikzpicture}
	% the origin
	\coordinate (O) at (0,0);
	\draw[-] (-3.5,0) -- (-1,0) node[right] {};
	\draw[->] (-0.5,0) -- (3.5,0) node[right] {$\hbox{Re}$};
	\draw[->] (0,-3.5) -- (0,3.5) node[above] {$\hbox{Im}$};
	% the circle and the dot at the origin
	\draw (O) node[circle,inner sep=1.5pt,fill] {} circle [radius=3cm];
	% the ``\theta'' arc
	\foreach \Point in {(120.5:2.99cm), (200.5:2.99cm), (339.5:2.99cm)}{
		\node at \Point {\textbullet};
	}
	\foreach \Point in {(160.5:2.99cm),(19.5:2.99cm)}{
	\node at \Point {o};}
	\node (b1) at (340:3cm) {};
\node (b2) at (160:3cm) {}; 
\draw[thick,dotted] (b1) to node {} (O);
	\draw 
	(200:3cm) coordinate (xcoord) -- 
	node[midway,below] {$r$} (O) -- 
	(120:3cm) coordinate (slcoord)
	pic [draw,<-,angle radius=1.7cm,"$\frac{2\pi}{\beta}$",above] {angle = slcoord--O--xcoord};
	\draw 
    (200:3cm) coordinate (xcoord1) -- 
    node[midway,below] {} (O) -- 
    (160:3cm) coordinate (slcoord1)
    pic [draw,<-,angle radius=1.3cm,"$\frac{\pi}{\beta}$"] {angle = slcoord1--O--xcoord1};	
	\node (a1) at (200:3cm) {};
	\node (a2) at (160:3cm) {}; 
	\draw[thick] (a1) to node {} (a2);	
	% the outer ``s'' arc
	\draw[|->]
	(-3cm-10pt,0)
	arc[start angle=180,end angle=200,radius=3cm+10pt]
	node[midway,fill=white] {$\sigma$};
		\draw[|->]
	(200:3cm+10pt)
	arc[start angle=200,end angle=270,radius=3cm+10pt]
	node[midway,fill=white] {$\angle M$};
	\end{tikzpicture}
	\caption{When $M\in\mathcal{M_{\text{odd}}}$, the limitation is independent of the period as when $T=\beta$, the constraint is activated by the opposite element. As a result, the limit for $\sigma$ is achieved when $2\sigma=\pi/\beta$ regardless of the period $T$.}\label{fig:2}
\end{figure}

\begin{remark}
	The last part of the statement is trivially true from the definition of the multipliers.
\end{remark}
 
\begin{IEEEproof}
		Let us considered the case $M\in\mathcal{M}$. For the sake of simplicity, let $G(e^{j\omega_1})=re^{j(-\pi+\sigma)}$ for $r>0$ and $0<\sigma\le\pi/2$. Then the limitation~\eqref{eq:one_f_w} holds for all $\sigma$ in the interval $[0,\pi/T]$, where the period $T$ of the complex exponential sequence is given in Lemma~\ref{lemma:periodic}. The limit case is shown in Fig.~\ref{fig:1} when $2\sigma=2\pi/T$. A maximum allowed phase for $G(e^{j\omega})$ to show that there is no suitable Zames-Falb multiplier implies in turn a maximum phase for Zames-Falb multipliers given by $\pi/2-\sigma$. As two different periods for $T$ are given in Lemma~\ref{lemma:periodic}, then the positive limits of both conditions \eqref{eq:9a} and \eqref{eq:9b} are found. 
		
		The same argument can be used when $G(e^{j\omega_1})=re^{j(\pi-\sigma)}$ for $r>0$ and $0<\sigma\le\pi/2$; see the limit case $2\sigma=2\pi/T$ in Fig~\ref{fig:3}, where $T$ is again defined by Lemma~\ref{lemma:periodic}. We find the minimum allowed phase for $G(e^{j\omega_1})$ which in turn can be translated into the negative limits of both conditions \eqref{eq:9a} and \eqref{eq:9b}. 
		 
		The same approach is used when $M\in\mathcal{M_{\text{odd}}}$. The only important difference in this case is that the periodicity of $e^{-j\omega_{1}i}$ is no longer relevant as condition~\eqref{eq:one_f_odd_w} include both positive and negative signs. Then the period of the exponential sequence is $2\beta$ regardless of $\alpha$. 
		
		The case $\alpha$ even is depicted in Fig.~\ref{fig:2}.
\end{IEEEproof}

	%The result follows from Corollary~\ref{cor:single} when it is applied to the rational frequency $\omega_1=\frac{\alpha}{\beta}$.  Let $\sigma$ be the phase of $G(e^{j\omega_{1}}-1/k)$, and let us assume that $\sigma\neq 0$ and $\sigma\in[-\frac{\pi}{2},\frac{\pi}{2}]$ without loss of generality.
	
	%Then, the set of complex numbers in the right side of~\ref{eq:one_f_odd_w} is given by $beta$ (even $\alpha$) or $2\beta$ (odd $\alpha$) points in the circle of radius $|G(e^{j\omega_{1}}-1/k)|$ with phases $\sigma+\frac{i}{\beta}\pi$ with $i=1,2,\dots,2\beta$ (odd $\alpha$) or $\sigma+\frac{2i}{\beta}$ (even $\alpha$) for $i=i=1,2,\dots,\beta$.

The obtained phase limitations are illustrated in Figs.~\ref{fig:pl} and~\ref{fig:pl_odd}.
\subsection{Tightness of the phase limitation}

Here we show that the conditions obtained in Theorem~\ref{Th:PL} are tight; at each frequency point it is possible to construct a single parameter Zames-Falb multiplier that meets the phase limitations in Theorem~\ref{Th:PL} with equality.

\begin{figure}[h]
	\centering
	\begin{tikzpicture}
	
\node (a) at (0,-1) {{\large $\frac{0}{1}$}};
\node (b) at (8,-1) {{\large$\frac{1}{1}$}};
\node (c) at (4,-2.5) {{\large$\frac{1}{2}$}};
\node (d) at (0,-2.5) {{\large$\frac{0}{1}$}};
\node (e) at (8,-2.5) {{\large$\frac{1}{1}$}};
\node (f) at (4,-4) {{\large$\frac{1}{2}$}};
\node (g) at (0,-4) {{\large$\frac{0}{1}$}};
\node (h) at (8,-4) {{\large$\frac{1}{1}$}};
\node (i) at (2,-4) {{\large$\frac{1}{3}$}};
\node (j) at (6,-4) {{\large$\frac{2}{3}$}};
\node (k) at (0,-5.5) {{\large$\frac{0}{1}$}};
\node (l) at (8,-5.5) {{\large$\frac{1}{1}$}};
\node (m) at (2,-5.5) {{\large$\frac{1}{3}$}};
\node (n) at (6,-5.5) {{\large$\frac{2}{3}$}};
\node (o) at (1,-5.5) {{\large$\frac{1}{4}$}};
\node (p) at (7,-5.5) {{\large$\frac{3}{4}$}};
\node (q) at (3,-5.5) {{\large$\frac{2}{5}$}};
\node (r) at (5,-5.5) {{\large$\frac{3}{5}$}};
\node (s) at (4,-5.5) {{\large$\frac{1}{2}$}};
\node (k1) at (0,-7) {{\large$\frac{0}{1}$}};
\node (l1) at (8,-7) {{\large$\frac{1}{1}$}};
\node (m1) at (2,-7) {{\large$\frac{1}{3}$}};
\node (n1) at (6,-7) {{\large$\frac{2}{3}$}};
\node (o1) at (1,-7) {{\large$\frac{1}{4}$}};
\node (p1) at (7,-7) {{\large$\frac{3}{4}$}};
\node (q1) at (3,-7) {{\large$\frac{2}{5}$}};
\node (r1) at (5,-7) {{\large$\frac{3}{5}$}};
\node (s1) at (4,-7) {{\large$\frac{1}{2}$}};
\node (aa) at (0.5,-7) {{\large$\frac{1}{5}$}};
\node (bb) at (1.5,-7) {{\large$\frac{2}{7}$}};
\node (cc) at (2.5,-7) {{\large$\frac{3}{8}$}};
\node (dd) at (3.5,-7) {{\large$\frac{3}{7}$}};
\node (ee) at (4.5,-7) {{\large$\frac{4}{7}$}};
\node (ff) at (5.5,-7) {{\large$\frac{5}{8}$}};
\node (gg) at (6.5,-7) {{\large$\frac{5}{7}$}};
\node (hh) at (7.5,-7) {{\large$\frac{4}{5}$}};
\draw [thick] (a) to node {} (c)
(b) to node {} (c)
(d) to node {} (i)
(g) to node {} (o)
(c) to node {} (i)
(i) to node {} (q)
(i) to node {} (o)
(f) to node {} (q)
(f) to node {} (r)
(c) to node {} (j)
(e) to node {} (j)
(j) to node {} (r)
(j) to node {} (p)
(h) to node {} (p)
(k) to node {} (aa)
(o) to node {} (aa)
(o) to node {} (bb)
(m) to node {} (bb)
(m) to node {} (cc)
(q) to node {} (cc)
(q) to node {} (dd)
(s) to node {} (dd)
(s) to node {} (ee)
(r) to node {} (ee)
(r) to node {} (ff)
(n) to node {} (ff)
(n) to node {} (gg)
(p) to node {} (gg)
(p) to node {} (hh)
(l) to node {} (hh)
;
\draw [thick, dotted] (a) to node {} (d)
(b) to node {} (e)
(d) to node {} (g)
(e) to node {} (h)
(c) to node {} (f)
(g) to node {} (k)
(i) to node {} (m)
(f) to node {} (s)
(j) to node {} (n)
(h) to node {} (l)
(k) to node {} (k1)
(l) to node {} (l1)
(m) to node {} (m1)
(n) to node {} (n1)
(o) to node {} (o1)
(p) to node {} (p1)
(q) to node {} (q1)
(r) to node {} (r1)
(s) to node {} (s1);
%(b) to node {2} node [swap] {2’} (c)
%(c) to node {3} node [swap] {3’} (a);
	\end{tikzpicture}

	\caption{First five sequences of the left branch of the Stern-Brocot tree. For each pair of neighbours at some level, a new element is generated by the mediant between to neighbours for the next level, i.e. the neighbours $p_k/q_k$ and $p_{k+1}/q_{k+1}$ generates a new element in the level below given by $(p_k+p_{k+1})/(q_k+q_{k+1})$. Adapted from~\cite{Concrete} and~\cite{wiki:2020}. }

	\label{Stern-Brocot}

\end{figure}
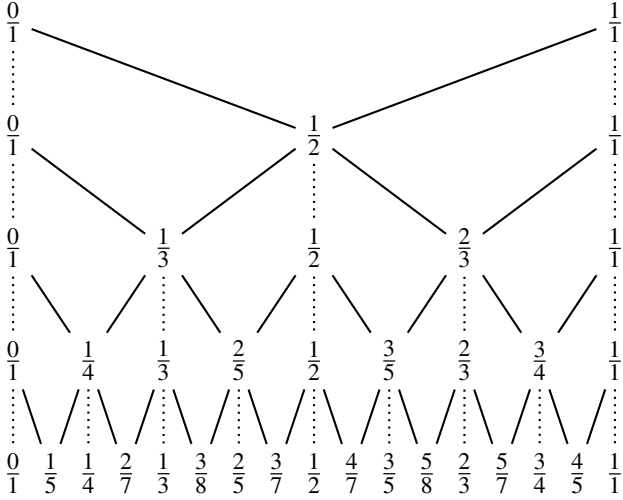

The result uses the left branch of the Stern-Brocot tree~\cite{Concrete} (see Fig.~\ref{Stern-Brocot}) so we can construct the required multiplier to achieve an arbitrary constraint. There is a strong relationship between the Stern-Brocot tree and Euclid's algorithm to find the B\'{e}zout coefficient of coprime numbers~\cite{Concrete}, and some of the results of this section are not restricted to the use of the Stern-Brocot tree, but it allows us to provide a closed form expression of the required multipliers. Furthermore, it allows us to construct a multiplier arbitrarily close to $\pm\pi/2$ at frequencies of the form $\gamma\pi$ with $\gamma$ irrational.

Loosely speaking, the sequence at a given level of the Stern-Brocot tree is denoted by $\{p_i/q_i\}_{i=1}^N$, with some abuse of notation. For instance, the first level is given by $\{0/1,1/1\}$, the second level is $\{0/1,1/2,1/1\}$, etc. 

We will require the following properties of the Stern-Brocot tree.

\begin{lemma}[B\'{e}zout identity~\cite{Concrete}] Any fraction $a/b$ with $a\perp b$ belongs to sequences of the Stern-Brocot tree. Moreover, for any two consecutive elements of a given sequence , $p_k/q_k$ and $p_{k+1}/q_{k+1}$,   
then 
\begin{equation}\label{eq:Bezout}
p_{k+1}q_k-p_kq_{k+1}=1.
\end{equation}	
\end{lemma}

\begin{lemma}\label{lemma:mult} Let us consider two consecutive elements of a given sequence for the tree, $p_k/q_k$ and $p_{k+1}/q_{k+1}$, with $p_k$ even. Then, the following properties hold:
	\begin{itemize}
		\item for the multiplier $M_1(z)=1+z^{-q_k}$,
		 $$\angle M_1(e^{\frac{p_{k+1}}{q_{k+1}}\pi j})=\frac{\pi}{2}\left(1-\frac{1}{q_{k+1}}\right),$$
		\item for the multiplier $M_2(z)=1+z^{q_k}$,
		 $$\angle M_2(e^{\frac{p_{k+1}}{q_{k+1}}\pi j})=-\frac{\pi}{2}\left(1-\frac{1}{q_{k+1}}\right),$$
		\item for the multiplier $M_3(z)=1-z^{q_{k+1}}$,
		 $$\angle M_3(e^{\frac{p_{k}}{q_{k}}\pi j})=\frac{\pi}{2}\left(1-\frac{1}{q_{k}}\right),$$
		\item for the multiplier $M_4(z)=1-z^{-q_{k+1}}$,
        $$\angle M_4(e^{\frac{p_{k}}{q_{k}}\pi j})=-\frac{\pi}{2}\left(1-\frac{1}{q_{k}}\right),$$
        \item for the multiplier $M_5(z)=1-z^{2q_{k+1}}$,
        $$\angle M_5(e^{\frac{p_{k}}{q_{k}}\pi j})=-\frac{\pi}{2}\left(1-\frac{2}{q_{k}}\right),$$
        \item for the multiplier $M_6(z)=1-z^{-2q_{k+1}}$
        $$\angle M_6(e^{\frac{p_{k}}{q_{k}}\pi j})=-\frac{\pi}{2}\left(1-\frac{2}{q_{k}}\right).$$        
	\end{itemize}
\end{lemma}

\begin{figure}
	\centering
\begin{tikzpicture}[scale=5]
% \draw[step=.5cm, gray, very thin] (-0.1,-1.2) grid (1.2,0.2); 
%\filldraw[fill=green!20,draw=green!50!black] (0,0) -- (3mm,0mm) arc (0:30:3mm) -- cycle; 
\draw[->] (-0.1,0) -- (1.2,0) node[right] {$\omega$};
\draw[->] (0,-0.8) -- (0,0.1) node[above] {$\angle M$};
%\draw (0,0) circle (1cm);
%\draw[very thick,red] (30:1cm) -- node[left,fill=white] {$\sin \alpha$} (30:1cm |- x axis);
\draw[very thick,blue] (0.75,-0) circle[radius=0.05 em] (0.05,-0.70) -- (0.74,-.01) ;
\draw[very thick,red] (0.25,-0) circle[radius=0.05 em] (0.26,-0.01) -- (1,-0.99);
\draw[very thick,dotted] (0,-0.5)--(0.25,-0.5)--(0.25,0);
\draw[very thick,dotted] (0,-0.66)--(0.75,-0.66)--(0.75,0);
 \foreach \x/\xtext in {0.25/\frac{p_k}{q_k}\pi, 0.75/\frac{p_{k+1}}{q_{k+1}}\pi} 
\draw (\x cm,2pt) -- (\x cm,0.75pt) node[anchor=south,fill=white] {$\xtext$};
 \foreach \y/\ytext in {-0.5/\frac{\pi}{2}\left(1-\frac{1}{q_k}\right), -0.66/\frac{\pi}{2}\left(1-\frac{1}{q_{k+1}}\right),0/\frac{\pi}{2}} 
\draw (1pt,\y cm) -- (-1pt,\y cm) node[anchor=east,fill=white] {$\ytext$};
\end{tikzpicture}
\caption{If $p_k$ is even, phase of the multipliers $1+z^{-q_{k}}$ (in red) and $1-z^{q_{k+1}}$ (in blue). If $p_{k}$ is even, phase of the multipliers $1-z^{-q_{k}}$ (in red) and $1+z^{q_{k+1}}$ (in blue). The multiplier in blue and the multiplier in red reach the phase limitation at frequencies $\frac{p_k}{q_k}\pi$ and $\frac{p_{k+1}}{q_{k+1}}\pi$, respectively.}\label{fig:mult}
\end{figure}
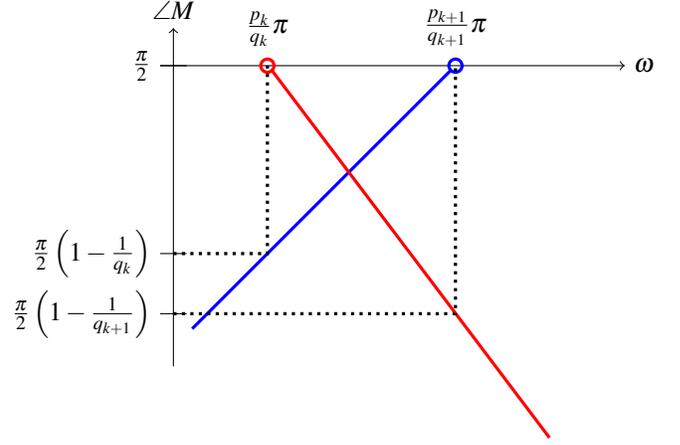

\begin{IEEEproof} The proof is based on B\'{e}zout's identity. For example, let us consider the multiplier $M_1$. Its phase is undefined at $\pi p_k/q_k$, and its slope is $-q_k/2$ when it is defined (see Fig.~\ref{fig:mult}). As a result, it is straightforward that
	$$\angle M_1(e^{\omega j})=\frac{\pi}{2}-(\omega-\pi \frac{p_k}{q_k})\frac{q_k}{2}, \quad \forall  \omega\in\left(\frac{p_k}{q_k}\pi,\frac{p_k+2}{q_k}\pi\right).$$
	
Substituting $\omega=\frac{p_{k+1}}{q_{k+1}}\pi$, it yields
	$$\angle M_1(e^{\frac{p_{k+1}}{q_{k+1}}\pi j})=\frac{\pi}{2}\left(1-\left(\frac{p_{k+1}}{q_{k+1}}-\frac{p_{k}}{q_{k}}\right)q_k\right).$$
	
	Then by using B\'{e}zout's identity~\eqref{eq:Bezout}, the result is 
	\begin{multline*}
		\angle M_1(e^{\frac{p_{k+1}}{q_{k+1}}\pi j})=\frac{\pi}{2}\left(1-\left(\frac{p_{k+1}q_k-p_{k}q_{k+1}}{q_{k}q_{k+1}}\right)q_k\right)=\\
		\frac{\pi}{2}\left(1-\left(\frac{1}{q_{k+1}}\right)\right).
	\end{multline*}
	
	The rest of the results can be obtained following the same approach.
\end{IEEEproof}

Then we can generate a multiplier reaching the limitation as follows.	
\begin{theorem}\label{th:3} Let $\alpha\in\mathds{Z^+},\beta\in\mathds{Z^+}$ with $\alpha\perp \beta$, $\alpha$ even, and $\alpha<\beta$. There is at least one Zames--Falb multiplier $M\in\mathcal{M}$ with phase $\frac{\pi}{2}(1-\frac{2}{\beta})$ at frequency $\frac{\alpha}{\beta}\pi$. Similarly, there is at least one Zames-Falb multiplier $M\in\mathcal{M_{\text{odd}}}$ with phase $\pm\frac{\pi}{2}(1-\frac{1}{\beta})$ at frequency $\frac{\alpha}{\beta}\pi$.
\end{theorem}

\begin{IEEEproof}
	The result follows from finding the neighbours of $\alpha/\beta$ in the Stern-Brocot tree, i.e. finding their B\'{e}zout's coefficients of $\alpha$ and $\beta$, and then applying the multipliers in Lemma~\ref{lemma:mult} as follows: if $\alpha$ is even, then take $p_k=\alpha$ and $q_k=\beta$ and the multipliers $M_3$ and $M_4$; if $\alpha$ is odd, then take $p_{k+1}=\alpha$ and $q_{k+1}=\beta$ and the multipliers $M_1$ and $M_2$.  
\end{IEEEproof}	  

\begin{figure}[ht]
	\centering
	\includegraphics[width=\linewidth]{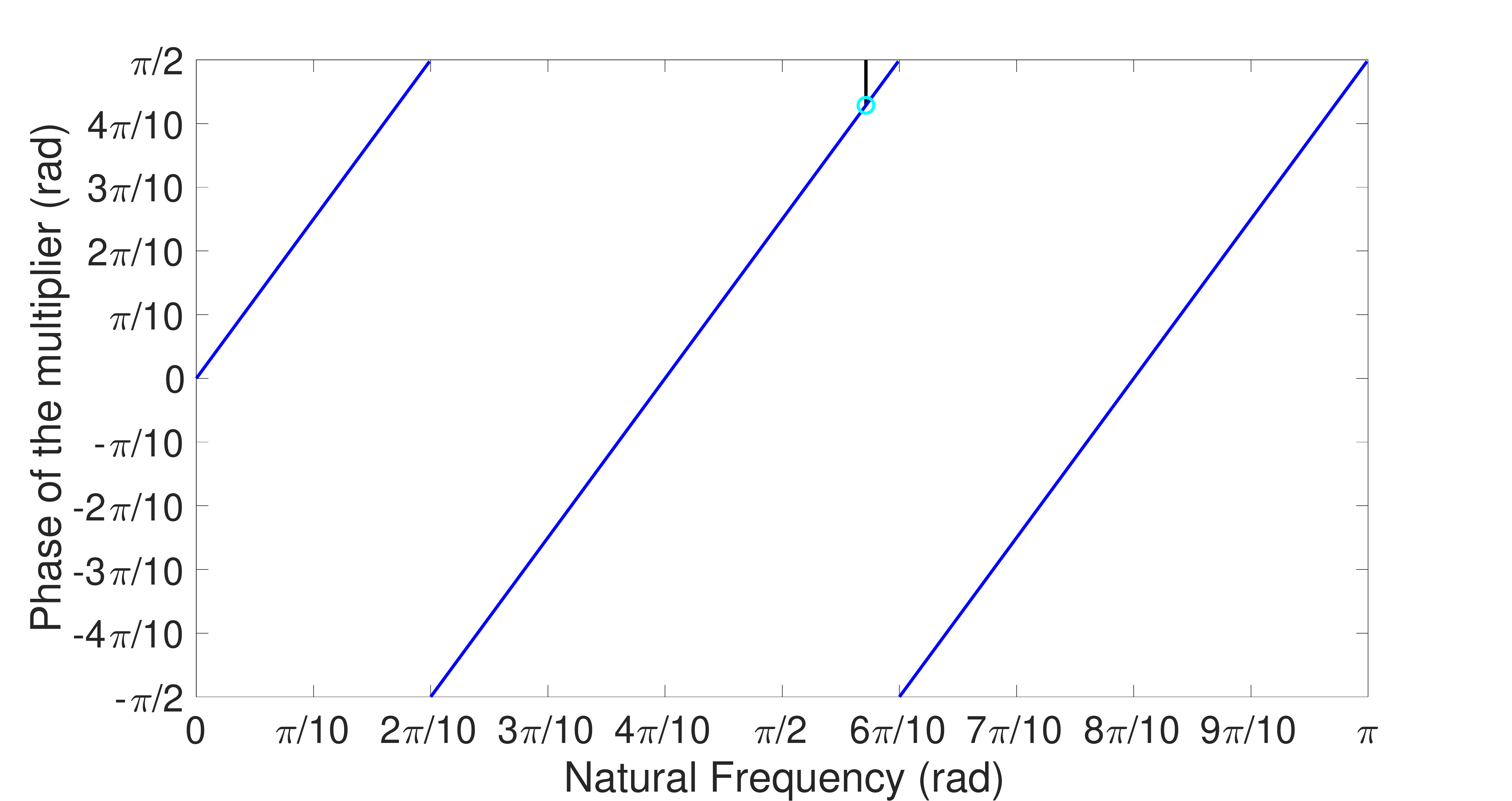}
	\caption{The multiplier $1+z^{5}$ reaches the phase $\frac{\pi}{2}\left(1-\frac{1}{7}\right)$ at the frequency $\frac{4\pi}{7}$. Note that the multiplier reaches infinite limitations, indeed they can be parametrised $\omega=\frac{p_k}{q_k}\pi$ where $p_k/q_k$ are all the neighbours of $\pi/5$, $3\pi/5$  for all the levels of the Stern-Brocot tree.}
	\label{fig:example}
\end{figure}

\paragraph{Example} Let us take $\alpha=4$ and $\beta=7$. It is straightforward to develop a bisectional search of the neighbours of $4/7$, so we can find that the right neightbour is $3/5$ in the first level of the tree when $4/7$ turns up. Then the multipliers $1+z^{\pm5}$ reach the phase $\pm\frac{\pi}{2}\left(1-\frac{1}{7}\right)$ at the frequency $\frac{4\pi}{7}$, see Fig.~\ref{fig:example}.
$\square$

Similarly, this development provides  another set of tools to show that the limitation collapses to $\pm\pi/2$ for irrational frequencies.

\begin{theorem} Let $\gamma\in \mathds{R^+}\backslash \mathds{Q}$, then for a given arbitrarily small $\epsilon>0$, there are  Zames-Falb multipliers of the form $M(z)=1-z^{\pm n_0}$ with phase larger than $\pi/2-\epsilon$ and smaller than  $-\pi/2+\epsilon$ at the frequency $\gamma\pi$.   
\end{theorem} 
\begin{IEEEproof}
	For a given $\epsilon$ take the Stern-Brocot tree with enough depth such that 
	$$\frac{p_k}{q_k}<\gamma<\frac{p_{k+1}}{q_{k+1}}$$
    with $\pi/(2q_k)<\epsilon$ and $\pi/(2q_{k+1})<\epsilon$. 
    
    Let us assume $p_{k+1}$ is even, then the phase of the multiplier $M_+(z)=1-z^{q_{k+1}}$ is in the interval $(\pi/2(1-1/q_k)),\pi/2)$ for any frequency between $\pi p_k/q_k$ and $\pi p_{k+1}/q_{k+1}$, hence the phase of the multiplier $M_+(z)=1-z^{q_{k+1}}$ is larger than $\pi/2-\epsilon$ at the frequency $\gamma\pi$. 
    Similarly, the multiplier $M_-(z)=1-z^{-q_{k+1}}$ is smaller than $-\pi/2+\epsilon$ at the frequency $\gamma\pi$ as $\angle M_-(z)=-\angle M_+(z)$ for all $|z|=1$. 
    
    If $p_{k+1}$ is odd, then $p_k$ is even, and the same argument can be followed for the multipliers $M_{\pm}=1-z^{\mp q_{k}}$.
\end{IEEEproof}

As a result, the single frequency limitation in Theorem~\ref{Th:PL} it is the best possible limitation for discrete-time Zames--Falb multipliers at a single frequency.

\subsection{Duality bounds with multiple frequencies}

As highlighted by J\"{o}nsson, one of the issues with continuous-time results for duality bounds is to obtain a suitable grid of frequencies. We have developed a powerful algorithm which does not require any numerical search. However, this method can be improved in some instances by using a grid of frequencies. In discrete-time, the selection of the frequency grid is not critical as the resulting conditions in Theorem~\ref{th:monotone} can be efficiently tested via a linear program.

\begin{proposition}\label{pro:linprog} For a stable $G$ and integer $\beta\geq2$, let us consider a set of frequencies $\omega_r=\frac{r}{\beta}\pi$ for $r=1,2,\cdots,\beta-1$. For $i=0,1,\cdots,2\beta-1$,
let us define 
	\begin{equation}\label{eq:defvar}
	\Lambda=\begin{bmatrix}
	\lambda_1\\
	\lambda_2\\
	\vdots\\
	\lambda_{\beta-1}
	\end{bmatrix},\quad\mathbf{v}_i^-=\begin{bmatrix}
	\hbox{Re}\{(1-e^{-j\omega_{1}i})G(e^{j\omega_1})\}\\
	\hbox{Re}\{(1-e^{-j\omega_{2}i})G(e^{j\omega_2})\}\\
	\vdots\\
	\hbox{Re}\{(1-e^{-j\omega_{2}i})G(e^{j\omega_{\beta-1}})\}
	\end{bmatrix},\end{equation}
	and
	\begin{equation}\label{eq:defvar1}
	\mathbf{v}_i^+=\begin{bmatrix}
	\hbox{Re}\{(1+e^{-j\omega_{1}i})G(e^{j\omega_1})\}\\
	\hbox{Re}\{(1+e^{-j\omega_{2}i})G(e^{j\omega_2})\}\\
	\vdots\\
	\hbox{Re}\{(1+e^{-j\omega_{\beta-1}i})G(e^{j\omega_{\beta-1}})\}
	\end{bmatrix}.\end{equation}

	%Assume there exists $\Lambda\succeq0$ given by~\eqref{eq:defvar} with $|\Lambda|>0$ such that
	Assume there exists a nonzero $\Lambda\succeq 0$ such that
	\begin{equation}\label{eq:11}
	\Lambda^\top \mathbf{v}_i^-\leq 0 \text { for all } i=0,1,\cdots,2\beta-1,
	\end{equation}
	then there is no Zames--Falb multiplier $M\in\mathcal{M}$ such that
	\begin{equation}
	\hbox{Re}\left\{M(e^{j\omega})G(e^{j\omega})\right\}>0, \quad \forall \omega\in[0,\pi]. 
	\end{equation}
	Similarly, assume there exists a nonzero $\Lambda\succeq 0$ such that
	\begin{equation}
	\Lambda^\top \mathbf{v}_i^-\leq 0, \text{ and } \Lambda^\top \mathbf{v}_i^+\leq 0, \text{ for all } i=0,1,\dots,2\beta-1,
	\end{equation}
	then there is no Zames--Falb multiplier $M\in\mathcal{M}_\text{odd}$ such that
	\begin{equation}
	\hbox{Re}\left\{M(e^{j\omega})G(e^{j\omega})\right\}>0, \quad \forall \omega\in[0,\pi]. 
	\end{equation}
\end{proposition}

\begin{IEEEproof}
	If there is a nonzero $\Lambda\succeq0$, such that 	
	\begin{equation}\
	\Lambda^\top \mathbf{v}_i^-\leq 0 \text { for all } i=0,1,\cdots,2\beta-1,
	\end{equation}
	then it is true that it holds for all $i\in\mathds{Z}$ as $\mathbf{v}_i=\mathbf{v}_{i\pm n2\beta}$. Hence the sequence $\{\lambda_1,\lambda_2,...,\lambda_{\beta-1}\}$ contains a nonzero element and
	\begin{equation}
\sum_{r=1}^{\beta-1} \hbox{Re}\left\{\lambda_r G(e^{j\omega_r})\right\}
\le 
\sum_{r=1}^{\beta-1} \hbox{Re}\left\{\lambda_r G(e^{j\omega_r})e^{-j\omega_r i}\right\},
	\end{equation}
	for all $i\in\mathds{Z}$. As a result, the condition \eqref{eq:1} in Thereom~\ref{th:monotone} are satisfied for all $M\in\mathcal{M}$.
	
	Similarly, if there exists a nonzero $\Lambda\succeq 0$ such that
	\begin{equation}
	\Lambda^\top \mathbf{v}_i^-\leq 0, \text{ and } \Lambda^\top \mathbf{v}_i^+\leq 0, \text{ for all } i=1,2,\dots,2\beta-1,
	\end{equation}
	then it is trivially true that it hold for all $i\in\mathds{Z}$ as $\mathbf{v}_i=\mathbf{v}_{i\pm n2\beta}$. Hence the sequence $\{\lambda_1,\lambda_2,...,\lambda_{\beta-1}\}$ contains a nonzero element and
\begin{equation}
\sum_{r=1}^{\beta-1} \hbox{Re}\left\{\lambda_r G(e^{j\omega_r})\right\}
\le 
\sum_{r=1}^{\beta-1} \hbox{Re}\left\{\lambda_r G(e^{j\omega_r})e^{-j\omega_r i}\right\},
\end{equation}
and
\begin{equation}
\sum_{r=1}^{\beta-1} \hbox{Re}\left\{\lambda_r G(e^{j\omega_r})\right\}
\le 
-\sum_{r=1}^{\beta-1} \hbox{Re}\left\{\lambda_r G(e^{j\omega_r})e^{-j\omega_r i}\right\},
\end{equation}
	for all $i\in\mathds{Z}$. As a result, the condition \eqref{eq:13} in Thereom~\ref{th:monotone} is satisfied for all $M\in\mathcal{M}_\text{odd}$.
\end{IEEEproof}

\section{Application to Absolute stability: Numerical results}\label{sec:results}

When absolute stability criteria are developed in the literature, the standard test consists of finding the maximum slope of the class of slope-restricted nonlinearities for a given stable system $G$~\cite{vidyasagar}. A loop transformation is normally used to obtain the following result:
\begin{corollary}\label{corollary:iqc_stability2}
	For the system in Fig. \ref{fig:lure}, let $G \in \mathbf{RH}_{\infty}$,  and $\phi$ is memoryless and $\phi\in S[0, k]$. Denote $\tilde{G}=G+\frac{1}{k}$. The system is $\ell_2$-stable, if there exists  $\Pi\in\Pi_{\phi}$ with $M\in\mathcal{M}$, such that 
	\begin{gather}\label{eq:iqc_stability2}
	\begin{bmatrix}
	\tilde{G}(e^{j\omega})  \\
	1 \\
	\end{bmatrix}^{\sim}
	\Pi(e^{j\omega})
	\begin{bmatrix}
	\tilde{G}(e^{j\omega})   \\
	1 \\
	\end{bmatrix}
	> 0,\quad  \forall\omega\in[0,\pi].
	\end{gather}
	Moreover, $M\in\mathcal{M_{\text{odd}}}$ in $\Pi$ if $\phi$ is also odd.
\end{corollary}

One of the first proposed solutions to the problem of absolute stability was proposed by Kalman~\cite{Kalman:1957} as follows:
\begin{definition}[Nyquist value, $k_N$] \label{df:Nyquist}
	The Nyquist value of a stable transfer function $G$ is
	\begin{equation*}
	k_N=\sup_{k}\{ k>0: (1+\tau kG(e^{j\omega}))^{-1} \ \textrm{is stable} \ \forall \tau\in[0,1] \}.
	\end{equation*}    
\end{definition}
\begin{conjecture}[Kalman conjecture~\cite{Kalman:1957}] Consider the feedback interconnection of $G$ and $\phi$. This feedback interconnection is asymptotically stable for all memoryless $\phi$ in $S[0,k]$ if and only if $k<k_N$.
\end{conjecture}

However, this solution has been proved wrong in general. For continuous-time, Fitts proposed the first counterexample~\cite{Fitts}; see~\cite{Leonov} for a discussion, and the conjecture has been proved to be true for first, second, and third order systems~\cite{Barabanov88}. In discrete-time, it is trivial that the conjecture is true for first order systems, but there are second-order counterexamples~\cite{Carrasco:2015,Heath:2015}.

\begin{figure}
	\centering
	\begin{picture}(300,270)
	\put(0,-75){\includegraphics[width=0.52\textwidth]{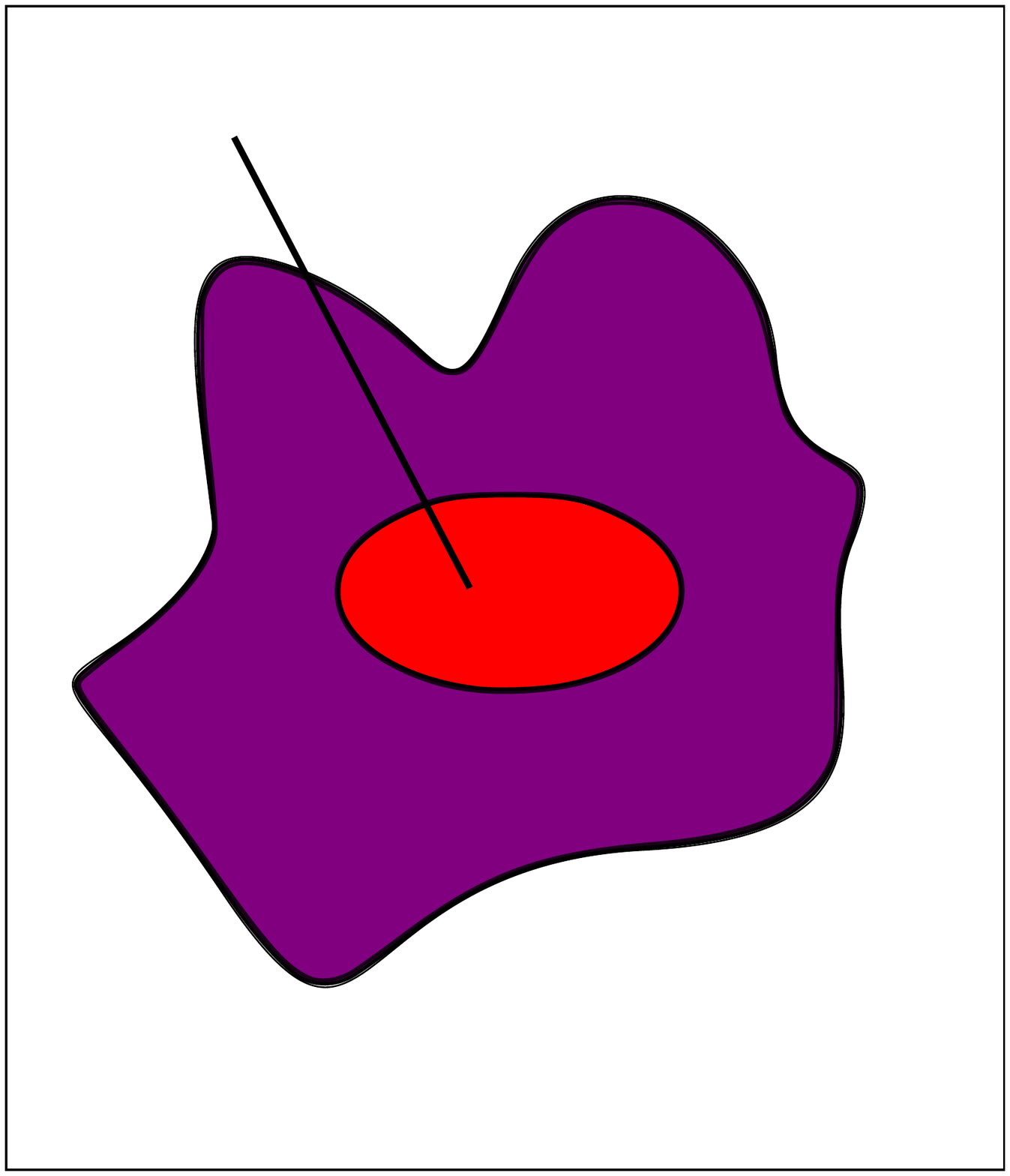}}
	\put(112,145){$G+\frac{1}{0^+}$}
	\put(60,245){$G$}
	\put(80,217){$G+\frac{1}{k_{ZF}}$}
	\put(25,45){$\mathbf{RH}_\infty$}
	\end{picture}
	\caption{Schematic representation of the trajectory of $G+1/k$ for stable plants with finite gain, i.e. $k_N<\infty$. In red, the set of LTI passive plants. In magenta, the set of LTI stable plants for which is possible to find a suitable Zames-Falb multiplier. It is worth highlighting that the gain at the crossing between the red and magenta areas is given by the Circle Criterion.}
	\label{fig:4}
\end{figure}

For a given $G$ with $k_N<\infty$, the set of plants $G+1/k$ is depicted in Fig.~\ref{fig:4}. A segment of this line lies within the subset of passive plants for $k\sim0$, i.e. the Zames-Falb multiplier $M(z)=1$ is suitable for these plants. A second segment is outside the set of passive plants but it is still possible to find a Zames-Falb multiplier. As $k_N<\infty$, there is no suitable Zames-Falb multiplier for $G$.

Mathematically, we can define the following critical gains:
\begin{itemize}
	\item supremum of the set of gains for which there exists a suitable LTI Zames-Falb multiplier, $k_{LTIZF}$;
	\item supremum of the set of gains for which there exists a suitable Zames-Falb multiplier, $k_{ZF}$; 
	\item supremum of the set of gains for which the system is absolutely stable, $k_{AS}$;
	\item supremum of the set of gains for which the system is stable in feedback interconnection with the linear gain, i.e. Nyquist gain, $k_N$.
\end{itemize}  

%\begin{figure}
%	\centering
%\begin{tikzpicture}
%\path 
%(0,0) rectangle (8,6) [draw]
%(0.5,5.5) node {$\mathbf{RH}_\infty$}
%(2.5,3) coordinate (A) node[above] {$\mathcal{A}$} ellipse (2 and 1.5)  [draw] +(0.5,-0.75) coordinate (A) node[below left] {$G+\frac{1}{0}$} -- (6,4.5) node[right]{{$G+\frac{1}{\infty}$}};
%\end{tikzpicture}
%\caption{Let $\mathcal{A}$ be the set of stable plants for which is possible to find a suitable Zames-Falb multiplier, then we ensure that the path crosses the boundary if $k_N<\infty$ for $G$.}
%\label{fig:4}
%\end{figure}

Trivially, $k_{LTIZF}\leq k_{ZF}\leq k_{AS}\leq k_{N}$\footnote{In the literature, two parallel problems are normally considered, when $\phi$ is slope-restricted and when $\phi$ is slope-restricted and odd, hence two sets of the above constants can be defined.}.
Then one can think of the search algorithm in~\cite{Joaquin:2019} as a lower bound of $k_{LTIZF}$, in short $\underline{k}_{LTIZF}$. The tools developed in this paper allow us to provide an upper bound of $k_{LTIZF}$, in short $\bar{k}_{LTIZF}$. For discrete-time examples, we shall show that $\bar{k}_{LTIZF}-\underline{k}_{LTIZF}$ is negligible for all tested examples; hence the current state of the art provides a very good estimation of $k_{LTIZF}$.

The results in~\cite{Joaquin:2019} show the best available lower bound of the solution of the problem, i.e. $k_{AS}$. We repeat here the conjecture proposed in~\cite{Shuai:2018}:

	\begin{conjecture}\label{cj:conjecture1}
		For an LTI $G$ and $0<k<k_N$, if there is no suitable LTI multipliers for $G+1/k$, then the Lurye system is not absolutely stable, i.e.
		$$k_{LTIZF}=k_{ZF}=k_{AS}.$$
	\end{conjecture}

This conjecture would imply that the set of linear time-variant (LTV) Zames--Falb multipliers can be ``phase-substituted'' by LTI Zames--Falb multipliers, hence it is not required to develop new techniques to find LTV Zames-Falb multipliers. On the other hand, the necessity of the existence of an LTI multiplier for absolute stability would imply there is no need for new stability criterion.

When the nonlinearity is slope restricted in the interval $S[0,k]$, then we apply a classical loop transformation as follows:
\begin{corollary} \label{th:sloperes}
	Let $G\in\mathbf{RH}_\infty$ and $k>0$. Assume there exists $0 < \omega_1\le \cdots < \omega_N \le \pi$, and $\lambda_1, \cdots, \lambda_N \ge 0$, where at least one $\lambda_r$ is nonzero. If
	\begin{multline}\label{eq:dual_con_nonodd_discrete}
	\sum_{r=1}^{N} \hbox{Re}\left\{\lambda_r\left(G(e^{j\omega_r})+\frac{1}{k}\right)\right\}\\
	\le 
	\min_{i\in \mathds{Z}}\left[\sum_{r=1}^{N} \hbox{Re}\left\{\lambda_r \left(G(e^{j\omega_r})+\frac{1}{k}\right)e^{-j\omega_r i}\right\}\right],
	\end{multline}
	then there is no Zames-Falb multipliers $M\in\mathcal{M}$ such that
	\begin{equation}
	\hbox{Re}\left\{M(e^{j\omega})\left(G(e^{j\omega})+\frac{1}{k}\right)\right\}>0, \quad \forall \omega\in[0,\pi]. 
	\end{equation}	
	Similarly, if 
	\begin{multline}\label{eq:dual_con_odd_discrete}
	\sum_{r=1}^{N} \hbox{Re}\left\{\lambda_r\left(G(e^{j\omega_r})+\frac{1}{k}\right)\right\}\\
	\le 
	-\max_{i\in \mathds{Z}}\left|\sum_{r=1}^{N} \hbox{Re}\left\{\lambda_r\left(G(e^{j\omega_r})+\frac{1}{k}\right)e^{-j\omega_r i}\right\}\right|,
	\end{multline}
	then there is no Zames-Falb multiplier $M\in\mathcal{M}_{odd}$ such that
	\begin{equation}
	\hbox{Re}\left\{M(e^{j\omega})\left(G(e^{j\omega})+\frac{1}{k}\right)\right\}>0, \quad \forall \omega\in[0,\pi]. 
	\end{equation} 
\end{corollary}

\begin{IEEEproof}
	The proof is straightforward by classical loop transformation as in Corollary 1.
\end{IEEEproof}

\subsection{Dual bounds $k_{LTIZF}$ by single frequency}

Although the limitations only hold for a countable set of frequencies, the results are very powerful as they allow us to provide an upper bound of $k_{LTIZF}$ in closed-form.

\begin{definition} 	Let $G\in\mathbf{RH}_{\infty}$.  Let $\psi:\mathds{R}\times\mathds{N^+}\rightarrow\mathds{R}$ be the function
\begin{equation}
\psi(\omega,\beta)=\frac{-\tan\left(\frac{1}{\beta}\pi\right)}{R(\omega)\tan\left(\frac{1}{\beta}\pi\right)+I(\omega)} 
\end{equation} where	
	$$R(\omega)=\hbox{Re}\left\{G(e^{j\omega})\right\},\quad I(\omega)=\left|\hbox{Im}\left\{G(e^{j\omega})\right\}\right|.$$

\end{definition}

\begin{proposition}\label{prop:one_frequency_condition}
	Let $G\in\mathbf{RH}_{\infty}$, given $\omega_{1}=\frac{\alpha}{\beta}\pi$, where  $\alpha\in\mathds{Z^+},\beta\in\mathds{Z^+}$  with $\alpha\perp \beta$ and $\alpha<\beta$. If
	\begin{equation}\label{eq:one_f_condition}
	k=\begin{cases}
	\psi(\omega_1,2\beta) \quad \text{ when } \alpha\text{ is odd}\\
	\psi(\omega_1,\beta) \quad \text{ when } \alpha\text{ is even}
	\end{cases}
	\end{equation}	
	is positive, then there is no $M\in\mathcal{M}$ such that $$\hbox{Re}\left\{M(e^{j\omega})\left(G(e^{j\omega})+\frac{1}{k}\right)\right\}>0, \quad \forall \omega\in[0,\pi].$$  
	
	Similarly, if
	\begin{equation}\label{eq:one_f_odd_condition}
	k=\psi(\omega_1,2\beta)
	\end{equation}
	is positive, then there is no $M\in\mathcal{M}_{\text{odd}}$ such that $$\hbox{Re}\left\{M(e^{j\omega})\left(G(e^{j\omega})+\frac{1}{k}\right)\right\}>0, \quad \forall \omega\in[0,\pi].$$  
\end{proposition}

\begin{IEEEproof}
	The proof follows from the phase limitation of Zames-Falb multipliers in Corollary~\ref{th:sloperes}.
\end{IEEEproof}

\subsection{Examples}\label{sec:results_example}
Since the value of $k_{AS}$ is unknown, it is hard to judge the conservativeness of $k_{LTIZF}$ or $\bar{k}_{LTIZF}$ individually. Therefore, both $k_{LTIZF}$ and $\bar{k}_{LTIZF}$ are needed to check the conservativeness by the dual gap. In this section, the examples in Table~\ref{Tab:examples} are considered, for which the sufficiency of the Kalman conjecture is  wrong. In all results below, the superscripts "no" and "odd" indicate the cases with non-odd and odd nonlinearities respectively.  
\begin{table}[h]\centering
	\begin{tabular}{|l|c|}
		\hline 
		
		Plant  & $k_N$ \\ \hline  
		
		$G_1(z)=\frac{0.1z}{z^2-1.8z+0.81}$ & $36.10000$  \\  
		
		$G_2(z)=\frac{ z^4 - 1.5 z^3 + 0.5 z^2 - 0.5 z + 0.5}{4.4 z^5 - 8.957 z^4 + 9.893 z^3 - 5.671 z^2 + 2.207 z - 0.5}$ &  $7.90700$  \\  
		
		$G_3(z)=\frac{z^3 - 1.95 z^2 + 0.9 z + 0.05}{z^4-2.8z^3+3.5z^2-2.412z+0.7209}$ & $2.74550$ \\

		$G_4(z)=\frac{-2.265 z^4 - 2.428 z^3 - 0.2606 z^2 + 0.253 z + 0.04455}{ z^5 + 2.465 z^4 + 2.201 z^3 + 0.8429 z^2 + 0.1188 z + 0.0006787}$ & $1.23987$ \\

		$G_5(z)=\frac{-2.225z^5 + 3.239z^4 - 1.708z^3 + 0.517z^2 - 0.1603z + 0.03239}{z^6 - 1.825z^5 + 1.927z^4 - 1.226z^3 + 0.1525z^2 + 0.1836z - 0.05546}$ &  $0.51373$ \\

		$G_6(z)=\frac{-0.08658z + 0.007162}{z^2 + 1.415 z + 0.5523}$ & $37.36307$\\ \hline

	\end{tabular} 
	\caption{Examples with the Nyquist values.} \label{Tab:examples}
\end{table}

\subsection{Lower bounds of $k_{LTIZF}$}\label{sec:results_primal}
The maximum slopes are obtained by solving the primal optimisation with the FIR Zames-Falb multipliers and the factorisation method in~\cite{Shuai:2014,Joaquin:2019}. Here, we set the causal step and the anticausal step of FIR multipliers as $n_z$. The toolbox CVX~\cite{CVX} with SDP solver SDPT3~\cite{SDPT3} is used to solve Linear Matrix Inequalities (LMIs). All results in this section were calculated in Matlab R2018a with processor: Intel(R) Core(TM) i7-8700 CPU @ 3.20GHz and RAM: 16.0GB.  The less-conservative results are listed in Table \ref{Tab:results} with the corresponding order of multipliers. 

\begin{table}[h]\centering
	\begin{tabular}{|l|c|c|c|c|}
		\hline
		
		&  $\underline{k}_{\text{ LTIZF}}^\text{no}$  & $n_z$ & 	$\underline{k}_\text{ LTIZF}^\text{odd}$ & $n_z$  \\ \hline 
		
		Ex.1 & $13.028317$ & $6$ & $13.511322$ & $20$  \\ \hline 
		
		Ex.2 & $3.823996$ & $5$  & $3.824034$& $10$  \\ \hline 
		
		Ex.3 & $0.802714$ & $5$ & $1.105645$ & $2$ \\ \hline 
		
		Ex.4 & $0.846650$ & $5$  &  $0.987666$ & $2$ \\ \hline 
		
		Ex.5 & $0.374445$ & $10$ & $0.374484$ & $8$ \\ \hline 
		
		Ex.6 & $13.262027$ & $8$ & $22.686904$& $6$  \\ \hline

	\end{tabular} 
	\caption{Lower bounds of $k_\text{LTIZF}$.} \label{Tab:results}
\end{table}

%\subsection{Dual slope bounds}\label{sec:results_dual}

\subsection{Duality bounds $\bar{k}_{LTIZF}$} \label{sec:results_dual_one_f}

Initially, the dual bounds of the class of Zames-Falb multipliers are solved by Proposition \ref{prop:one_frequency_condition}. The less-conservative results with the corresponding frequencies are in Table \ref{Tab:results0}.

\begin{table}[h]\centering
	\begin{tabular}{|l|c|c|c|c|}
		\hline		
		$\;$ &  $\bar{k}_\text{LTIZF}^\text{no}$ & $\omega_1^\text{no}$ &  $\bar{k}_\text{LTIZF}^\text{odd}$ &  $\omega_1^\text{odd}$ \\ \hline 
		
		Ex.1 & 13.028374 & $\frac{2}{7}\pi$ & \textbf{13.575410} &  $\frac{1}{3}\pi$ \\ \hline  
		
		Ex.2 & 3.824040 & $\frac{1}{2}\pi$&   3.824040 & $\frac{1}{2}\pi$ \\ \hline 
		
		Ex.3 &  0.802745 & $\frac{2}{5}\pi$ & 1.105649 & $\frac{1}{2}\pi$ \\ \hline
		
		Ex.4 & 0.846657 & $\frac{2}{3}\pi$ & 0.987671 &  $\frac{1}{2}\pi$   \\ \hline 
		
		Ex.5 &   0.374491 &$\frac{1}{3}\pi$ &   0.374491 &  $\frac{1}{3}\pi$\\ \hline  
		
		Ex.6 &  13.262035 & $\frac{2}{3}\pi$ & 22.686907 &  $\frac{1}{2}\pi$ \\ \hline

	\end{tabular} 
	\caption{Upper bound of $k_\text{LTIZF}$ using Proposition~\ref{prop:one_frequency_condition} with the corresponding frequencies. In bold the only result with significant differences with Table~\ref{Tab:results}.} \label{Tab:results0}
\end{table}

The dual gaps between $\underline{k}_{LTIZF}$ in Table~\ref{Tab:results} and $\bar{k}_{LTIZF}$ in Table~\ref{Tab:results0} are all very small, except Ex.1 for odd nonlinearities, where the dual gap is $0.474\%$. 

\begin{figure}
	\centering
	\includegraphics[width=\linewidth]{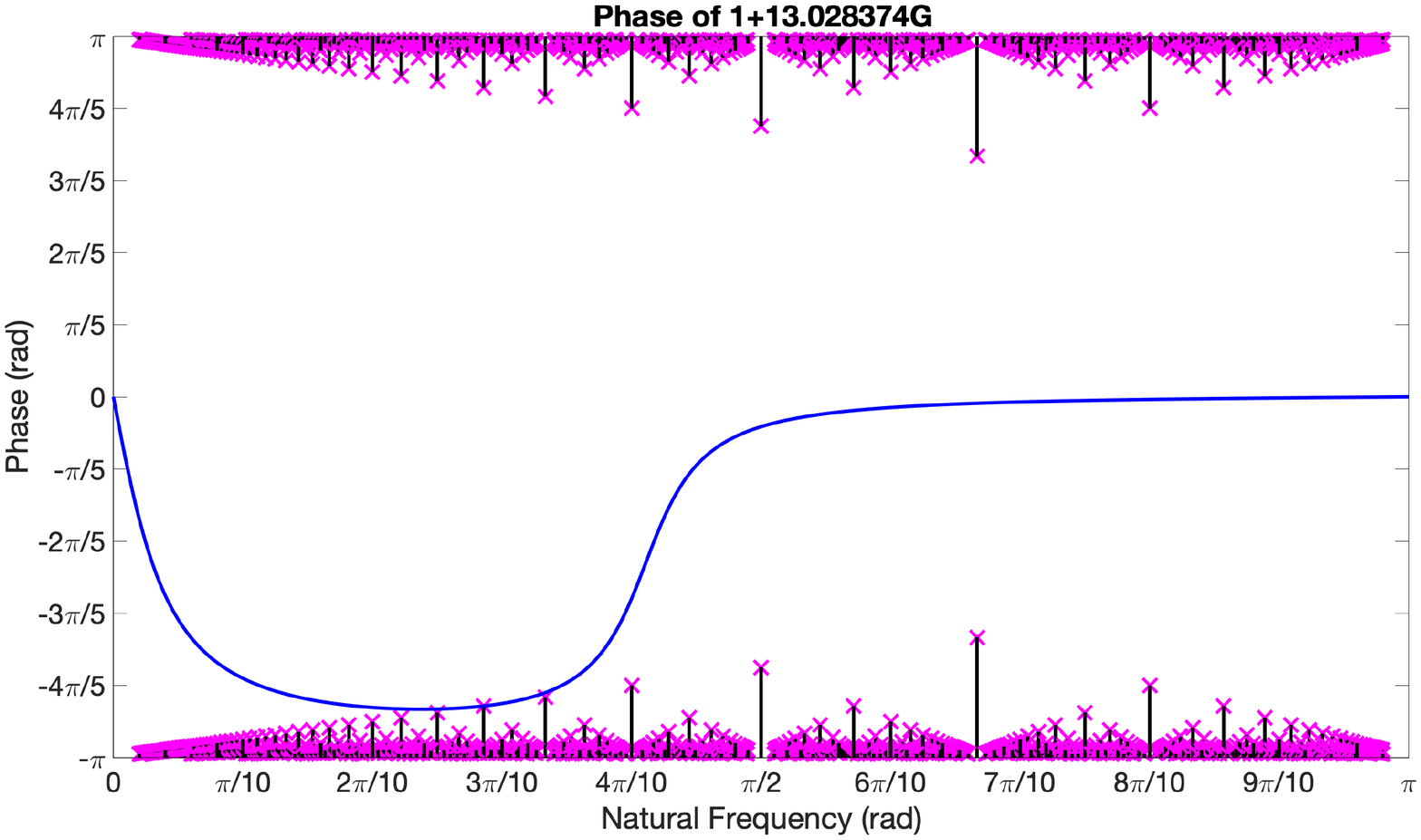}
	\caption{Example 1 reaching the limitations for non-odd nonlinearities when $k^\text{no}_\text{LTIZF}=13.028317$. The first limitation is reached at $\omega_{1}^\text{no}=\frac{2}{7}\pi$.}
	\label{fig:example1nonodd}
\end{figure}

\begin{figure}
	\centering
	\includegraphics[width=\linewidth]{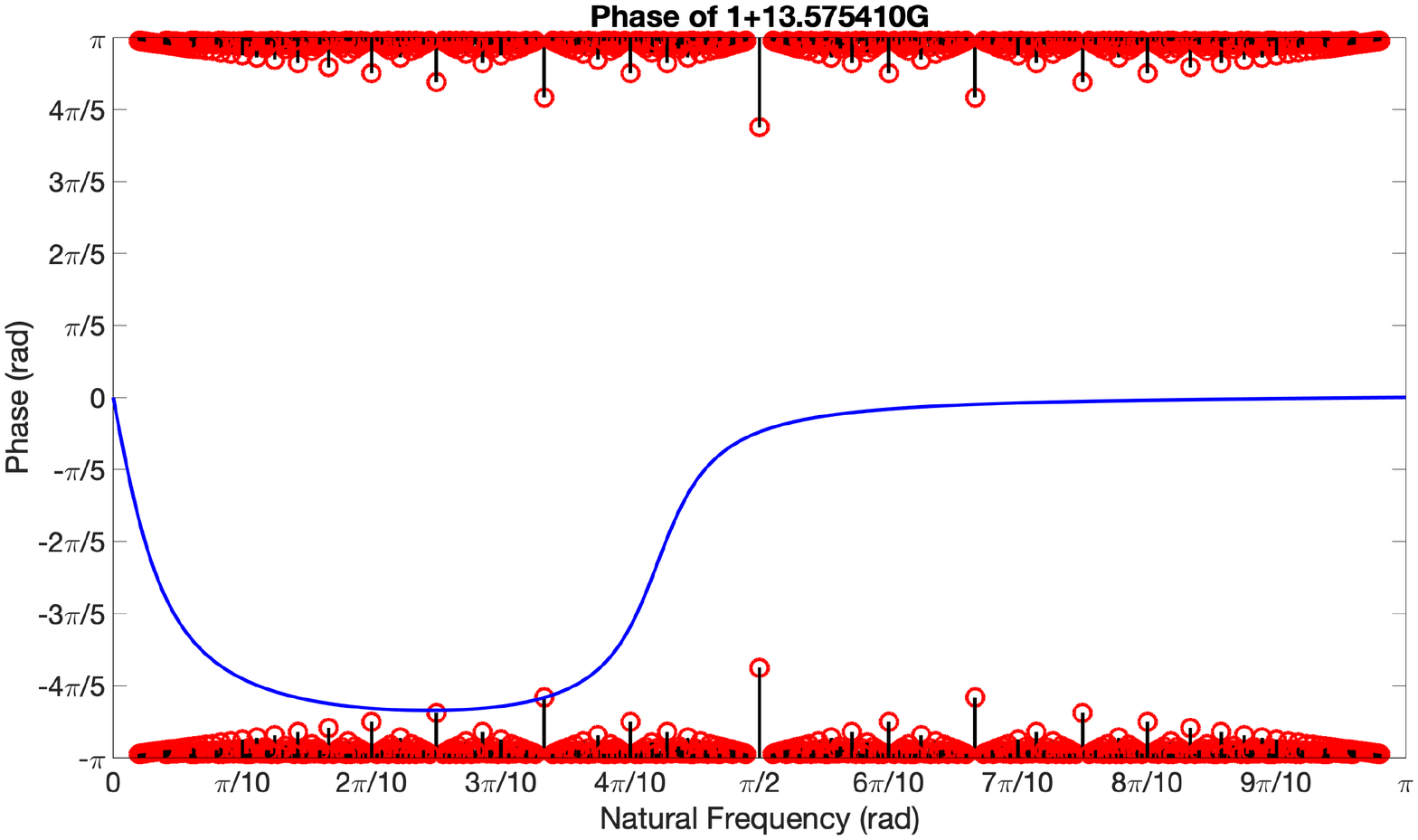}
	\caption{Example 1 reaching the limitations for odd nonlinearities when $k^\text{odd}_\text{LTIZF}=13.511322$. The first limitation is reached at $\omega_{1}^\text{odd}=\frac{1}{3}\pi$.}
	\label{fig:example1odd}
\end{figure}

It is possible to enhance the Bode plot with the limitations in Theorem~\ref{Th:PL} translated into the plant $1+kG$, so the user may explore graphically when the limitations are reached. For instance, if the available phase of the multiplier at $\omega=\pi/3$ is within the interval $(-\pi/3,\pi/3)$, then the Bode plot of $1+kG$ must lie in the interval $(-5\pi/6,5\pi/6)$ if $M(e^{j\omega})(1+kG(e^{j\omega}))>0$ at $\omega=\pi/3$.
As examples, the Bode plots of $(1+kG)$ at limiting values of $k$ for non-odd and odd cases in Table~\ref{Tab:results0}, are given in Figures~\ref{fig:example1nonodd} and~\ref{fig:example1odd}. These limitation have time-domain implications in terms of the existence of limit cycles~\cite{Seiler:2020}.

Finally, Proposition~\ref{pro:linprog} is used for Ex.1 for odd nonlinearities, with $\beta=250$. By using a standard bisectional algorithm, we obtain that  $\bar{k}^\text{odd}_{LTIZF}=13.511740$, reducing the dual gap to $0.003\%$. For the rest of the examples, Proposition~\ref{pro:linprog} does not improve Proposition~\ref{prop:one_frequency_condition}.

\subsection{Convexity properties}

A fundamental result  that can be obtained from the duality condition:

\begin{theorem}\label{th:2}
	The set of plants with a suitable Zames-Falb multiplier is not convex.
\end{theorem}

\begin{IEEEproof} We proof the statement by providing a numerical counterexample.

Consider the plant
$$G=0.2(G_1+1/12.9)+0.8(G_2+1/3.8).$$
From Table~\ref{Tab:results}, $G$ lies in the line connecting two plants with a suitable Zames-Falb multiplier. However, there is no suitable Zames-Falb for $G$ by applying Proposition~\ref{pro:linprog} with $\beta=40$.
\end{IEEEproof}

\subsection{Relations with the phase limitation approach}

%We compare our results with these given in~\cite{Shuai:2018} summarised in Section~\ref{sec:phase_limit}. When the results are tested for very similar frequencies, i.e. $a\to b$ in (\ref{eq:phase_limit}), we obtained
%\begin{equation*}
%\lim_{a\to b} \mu^{d}=\frac{|\sin(bn)|}{1-\cos(bn)},\quad 
%\lim_{a\to b} \mu^{d}_{odd}=\frac{|\sin(bn)|}{1-|\cos(bn)|},
%\end{equation*}
%Hence, although these two methods are closely related, they are different in the underpinning concepts.  

Theoretically, by involving a sufficiently large number of combinations of $0\le a<b\le \pi$, it is possible to obtain a least-conservative dual bound $\bar{k}_{LTIZF}$ by Theorem \ref{th:phase_limitation}.

% Then, it is natural to propose the conjecture below.
%\begin{conjecture}\label{cj:phase_limitation}
%	Assume that $\bar{k}_{D,min}$ is the least-conservative solution obtained by Proposition \ref{prop:one_frequency_condition}, where the corresponding frequency is $\omega$.  Then, the least-conservative solution $\bar{k}_{PL,min}$ obtained by Theorem \ref{th:phase_limitation} is more conservative than $\bar{k}_{D,min}$, where the corresponding frequencies $a$ and $b$  are close to $\omega$ on both sides respectively.
%\end{conjecture}

%In Section \ref{sec:results_phase_limitation}, this conjecture is true for all examples. 

The less-conservative dual bounds by Theorem \ref{th:phase_limitation} with Algorithm \ref{alg:phase_limit} (see Appendix~A) are listed in Table \ref{Tab:results2} together with the corresponding frequencies $a$ and $b$, where $a$ and $b$ are checked with the increment of $10^{-4}$ rad. Here, we choose the small resolution for less-conservative results at the cost of larger computational load.

\subsubsection{Conservativeness}
First of all, the numerical results here are less conservative than the algorithm~\cite{Shuai:2018} as there is no efficient algorithm. We have used a brute-force approach by using grid of frequencies. The algorithm is given in Appendix~\ref{Sec:appendix}.

\begin{table}[h]\centering
	\begin{tabular}{|l|c|c|c|c|c|c|}
		\hline
		
		$\;$ &  $\bar{k}^\text{no}_\text{LTIZF}$ &$a^\text{no}$& $b^\text{no}$& $\bar{k}_\text{LTIZF}^\text{odd}$ & $a^\text{odd}$& $b^\text{odd}$  \\ \hline 
		
		Ex.1 & $13.03978$ & $0.8966$ & $0.8986$ & $13.60146$  & $1.0455$ & $1.0489$\\ \hline  
		
		Ex.2 & $3.84105$ & $1.5684$ & $1.5732$ & $ 3.84820$  & $1.5674$ & $1.5742$\\ \hline  
		
		Ex.3 & $0.80372$ & $1.2532$ & $1.2601$ & $1.11120$  & $1.5674$ & $1.5742$\\ \hline  
		
		Ex.4 & $0.85026$ & $2.0902$ & $2.0986$ & $0.98853$  & $1.5674$ & $1.5742$\\ \hline  
		
		Ex.5 & $0.37586$ & $1.0460$ & $1.0484$ & $0.37644$  & $1.0455$ & $1.0489$\\ \hline  
		
		Ex.6 & $13.29765$ & $2.0902$ & $2.0986$ & $22.73067$  & $1.5674$ & $1.5742$\\ \hline

	\end{tabular} 
	\caption{Upper bounds computed with the phase limitation approach and the corresponding frequency pairs~\cite{Shuai:2018}.} \label{Tab:results2}
\end{table}

%
%\begin{table}[h]\centering
%	\begin{tabular}{|l|c|c|c|c|c|c|}
%		\hline
%		
%		$\;$ &  $\bar{k}^\text{no}_{PL,min}$ &$a^\text{no}$& $b^\text{no}$& $\bar{k}_{PL,min}^\text{odd}$ & $a^\text{odd}$& $b^\text{odd}$  \\ \hline 
%		
%		Ex.1 & $13.03912$ & $0.89666$ & $0.89854$ & $13.60130$  & $1.04551$ & $1.04889$\\ \hline  
%		
%		Ex.2 & $3.83401$  & $1.56938$ & $1.57221$  & $3.83815$  & $1.56880$ & $1.57279$   \\ \hline 
%		
%		Ex.3 & $0.80309$ & $1.25529$ & $1.25798$ & $1.10892$  & $1.56879$ & $1.57280$\\ \hline
%		
%		Ex.4 & $0.84717$ & $2.09379$ & $2.09500$ & $0.98818$  & $1.56878$ & $1.57281$ \\ \hline 
%		
%		Ex.5 & $0.37585$ & $1.04601$ & $1.04839$ & $0.37643$  & $1.04551$ & $1.04889$ \\ \hline  
%		
%		Ex.6 & $13.26713$ & $2.09379$ & $2.09500$ & $22.71274$  & $1.56879$ & $1.57280$ \\ \hline
%		
%		
%	\end{tabular} 
%	\caption{Less-conservative solutions by phase limitation approach} \label{Tab:results2}
%\end{table}

Next, comparing the results in Tables \ref{Tab:results0} and \ref{Tab:results2}, the results in~\cite{Shuai:2018} are more conservative than the one frequency results. In addition, $\omega_{1}\in(a,b)$ hence there is a clear connection between them, but a finite of between 2 and 5 mrad is required to extract the limitation~\cite{Shuai:2018}.

\subsubsection{Computational load}
Here, we compare the computational time of Ex. 1 with nonodd nonlinearities as an example. In phase limitation approach, the computational load depends exponentially on the initial bounds $k_m$ and $k_n$, the accuracy required on $k$, the number of $n$ to obtain the maximum in~(\ref{eq:phase_limit}), the resolution of $a$ and $b$, and the final value of $a$ and $b$. 

%For Ex.1, in Algorithm \ref{alg:phase_limit}, we initialise $k_n=0$ and $k_m=k_N$, and require $k_m-k_n\le10^{-5}$ in the end. We set $n=1:100$ in  (\ref{eq:phase_limit}). In addition, as the bisection search converges to a neighbour of the final value quickly, we can approximate there are  around $14000$ frequencies of the initial frequency grid that make $\hbox{Re}\left(G(e^{j\omega})+\frac{1}{k}\right)\leq0$ hold. Then, as $a<b$, there are at most $\frac{14000(14000-1)}{2}=97993000$ times of test on the phase limitation condition for each iteration of $k$. On the contrary, in Proposition \ref{prop:one_frequency_condition}, we consider the single frequency $\omega_1=\frac{1}{14000}\pi,\frac{2}{14000}\pi,\cdots,\pi$ individually, in order to include similar number of frequencies with the phase limitation approach.

For the parameters above, with a full search over the possible combination of frequencies, in Ex.1, it takes $16743.0447$  seconds  ($4$ hours $39$ minutes $3$ seconds) by the phase limitation approach. However, it takes 20 milliseconds to test 210 frequencies by Proposition \ref{prop:one_frequency_condition}. When Proposition~\ref{pro:linprog} is used, the linear program can be solved in $68.16$ seconds.

%Furthermore, as shown in Table \ref{Tab:results}, the denominators in the frequency are often $2,3,5,7$, which implies that $210$ individual frequencies from $\frac{1}{210}\pi$ to $\pi$ would be reliable to obtain a less-conservative result. In this way, it takes around $0.055$ seconds in each example by by Proposition \ref{prop:one_frequency_condition}. Finally, when the linear program is used, the solution is obtained in a few minutes. 

In short, compared with the phase limitation in~\cite{Shuai:2018}, the duality optimisation by Proposition \ref{prop:one_frequency_condition} is more competitive in terms of conservativeness with a huge difference in computational time.

\section{Conclusions}

The paper develops novel conditions to discard the existence of a suitable Zames-Falb multiplier. The result can be used at a single frequency, avoiding the need of an optimization; and at multiple frequencies by the use of a linear program. The application of both results to absolute stability generates an upper bound for the search in~\cite{Joaquin:2019}. In the tested examples, the larger gap between results is below $0.005\%$, and it is significant better than previous results in the literature. 

In conjunction with~\cite{Joaquin:2019}, the results in this paper allows us to conclude that the current state-of-the-art for discrete-time Zames-Falb multipliers algorithms are very efficient for both, the lower and upper bound for the existence of a Zames--Falb multiplier. However, other questions still deserve attention. For instance, linear-time varying (LTV) multipliers are included in the original definition by~\cite{Willems:1968}, but there are no numerical searches for them. It is still open whether LTV multipliers are useful or not. Moreover, the final aim of the problem is to find sufficient and necessary conditions for absolute stability. We conjecture that the lack of an LTI Zames-Falb multiplier implies that the system $G$ is not absolute stable.

Future research should focus now on the meaning of this value. We believe that Conjecture~\ref{cj:conjecture1} summaries the open questions for absolute stability of SISO systems. A partial answer has been provide in~\cite{Seiler:2020}, where the single frequency result provided in this paper has been connected to periodic behavior for the Lurye system. However, an interpretation of the multiple frequency condition is still open.

\section*{Acknowledgment}

The first author would like to thank the Department of Electrical and Electronic Engineering at The University of Manchester for its support.

This work was partially supported by EPSRC project EP/S03286X/1.

%\begin{figure}[ht]
%	\centering
%	\includegraphics[width=0.7\linewidth]{ex1_circle}
%\end{figure}
%
%
%\begin{figure}[ht]
%	\centering
%	\includegraphics[width=0.7\linewidth]{ex1_circle1}
%\end{figure}
%

%
%\begin{figure}[ht]
%	\centering
%	\includegraphics[width=0.7\linewidth]{ex1_circle2}
%\end{figure}
%

%\begin{figure}[ht]
%	\centering
%	\includegraphics[width=0.7\linewidth]{ex1_circle_two_fre}
%\end{figure}

%
%\bibliographystyle{IEEEtran}
%\bibliography{refer} 

\appendix \label{Sec:appendix}
 \subsection{Bisection search on $k$ for Theorem~\ref{th:phase_limitation}}\label{sec:appendix_A0}
This appendix describe the algorithm used to test phase limitations given in~\cite{Shuai:2018}.

\begin{algorithm}\label{alg:phase_limit} $\quad$
	\begin{enumerate}
		\item Preselect a frequency grid $0=\omega_1<\omega_2\cdots<\omega_N=\pi$.
		\item Initialise the upper and lower bounds of the slope as $k_m$ and $k_n$ respectively, 
		\item Update the slope $k=(k_m+k_n)/2$.
		\item Select the frequencies in Step $1$, such as $\omega_{d_1}$ to $\omega_{d_n}$, where $\hbox{Re}\left(G(e^{j\omega_{d_i}})+\frac{1}{k}\right)\leq0$, $\forall i=1\cdots n$. Denote the phase as $\sigma$.
		\item Obtain the ideal phase of the multiplier, $\eta$, as $\eta=-\frac{\pi}{2}-\sigma$ when $\sigma\le -\frac{\pi}{2}$, or $\eta=\frac{\pi}{2}-\sigma$ when $\sigma\ge \frac{\pi}{2}$.
		\item Set the scalar $a_i=\omega_{d_i}$ for every $i=1,\cdots, n-1$, and set the corresponding array $\pmb{b}_i=[
		\omega_{d_{i+1}}, \cdots,  \omega_{d_{n}}]$. 
		\item Calculate $\rho^d_{(odd)}$ at each $a=a_i$ combined with each $b$ in $\pmb{b}_i$. 
		If at a frequency pair $(a,b)$, $\rho^d_{(odd)} \le \eta$ $\forall \omega\in[a,b]$ when $\sigma\le -\frac{\pi}{2}$, or $\rho^d_{(odd)}\ge \eta$ $\forall \omega\in[a,b]$ when $\sigma\ge \frac{\pi}{2}$, then stop Step 7, and update $k_m=k$. 
		If the above condition cannot hold with any frequency pair $(a,b)$, then update $k_n=k$.	    
		\item Repeat from Step $3$ with the updated $k_m$ and $k_n$ until they are close enough, meanwhile the condition in Step 7 holds. The result is $\bar{k}_\text{LTIZF}=k_m$.% and the corresponding frequencies are $a^*$ and $b^*$.
	\end{enumerate}

\subsection{Continous-time duality results for Zames-Falb multipliers}

For completeness, we develop the duality condition for non-odd nonlinearities. For the contibuous-time domain notation of this Appendix, we refer to the reader to~\cite{Jonsson:1996,Joaquin:2016}.

The following result was developed by Jonss\"on and Laiou: 

\begin{theorem}[\cite{Jonsson:1996, Jonsson:1996c}]\label{th:duality_odd_continuous}
	Let $G\in\mathbf{RH}_\infty$. Assume there exist $0 < \omega_1< \cdots < \omega_{N-1} < \infty$, $\omega_N=\infty$ and non-negative $\lambda_1, \cdots, \lambda_N\ge 0$, where at least one $\lambda_r$ is nonzero, such that
	\begin{equation}\label{eq:dual_con_odd}
	\sum_{r=1}^{N} \hbox{Re}\left\{\lambda_r G({j\omega_r})\right\}
	\le -\sup_{t\in \mathds{R}}\left|\sum_{r=1}^{N-1}\hbox{Re}\left\{\lambda_r G({j\omega_r}) e^{-j\omega_r t}\right\}\right|,
	\end{equation}
	then there is no continuous-time Zames-Falb multipliers $M\in\mathcal{M}_{odd}$ such that  $\hbox{Re}\left\{M({j\omega})G({j\omega})\right\}>0$, $\forall \omega\in
	\mathds{R}\cup\{\pm \infty\}$. 
\end{theorem}

\begin{remark}
	We have adapted the original result to negative feedback configuration.  
\end{remark}

%	Following the similar analysis on $M\in \mathcal{M}$  for non-odd nonlinearities in Theorem \ref{th:lemma_4_1},  the continuous-time condition in Theorem \ref{th:duality_odd_continuous} can be relaxed to the non-odd case as follows.

In this result, the uncertainty $\phi$ in Lurye systems is assumed to be odd. We provide the duality bounds of continuous-time Zames-Falb multipliers, where the oddness of nonlinearities in Theorem~\ref{th:duality_odd_continuous} is relaxed.

\begin{theorem}\label{th:duality_nonodd_continuous}	
	Let $G\in\mathbf{RH}_\infty$. Assume there exists $0 < \omega_1\le \cdots \le \omega_{N-1} < \infty$, $\omega_N=\infty$ and $\lambda_1, \cdots, \lambda_N\ge 0$, where at least one $\lambda_r$ is nonzero, such that
	\begin{multline}\label{eq:dual_con_nonodd}
	\sum_{r=1}^{N} \hbox{Re}\left\{\lambda_rG({j\omega_r})\right\} 
	\le \inf_{t\in \mathds{R}}\left[\sum_{r=1}^{N-1} \hbox{Re}\left\{\lambda_r G({j\omega_r})e^{-j\omega_r t}\right\}\right],
	\end{multline}
	then there is no continuous-time Zames-Falb multiplier $M\in\mathcal{M}$ such that  $\hbox{Re}\left\{M({j\omega})G({j\omega})\right\}>0$, $\forall \omega\in
	\mathds{R}\cup\{\pm \infty\}$.
\end{theorem}

%\begin{remark} We have followed the same notation as J\"onsson. However, some reader may find useful to think of~\eqref{eq:dual_con_nonodd} as the following expression 
%\begin{equation}\label{eq:dual_con_nonodd_alter}
%\sum_{r=1}^{N-1} \hbox{Re}\left\{\lambda_rM_t(j\omega_r)G({j\omega_r})\right\} 
%\le0,
%\end{equation}
%where $M_t(j\omega)=1-e^{-j\omega t}$, for all $t>0$.	 
%\end{remark}

\begin{IEEEproof}
	The proof follows from Theorem 4.2 in~\cite{Jonsson:1996}. In the spirit of~\cite{Jonsson:1996c}, the duality condition is satisfied if  
	\begin{equation}\label{eq:th5-1}
	\sum_{r=1}^{N} \hbox{Re}\left\{\lambda_rG({j\omega_r})\right\} 
	\le \sum_{r=1}^{N-1} \hbox{Re}\left\{\lambda_rG({j\omega_r})e^{-j\omega_r t}\right\},
	\end{equation}
	for all $t\in\mathds{R}$. This condition implies that the left hand side must be non-positive. Then, for all signals $h\in\mathcal{L}_1(\mathds{R})$ such that $h(t)\ge 0$, and $\|h\|_1\le 1$, it follows
	\begin{equation}\label{eq:th5-2}
	\left(1-\int_{-\infty}^{\infty}h(t)dt\right)\hbox{Re}\left\{\lambda_rG({j\omega_r})\right\} 
	\le 0.
	\end{equation}
	
	Then, substituting (\ref{eq:th5-1}) to (\ref{eq:th5-2}) for each $t$ in the integration, it yields
	\begin{equation}
	\sum_{r=1}^{N} \hbox{Re}\left\{\lambda_rG({j\omega_r})\right\} 
	- \int_{-\infty}^{\infty}h(t)\sum_{r=1}^{N-1} \hbox{Re}\left\{\lambda_rG({j\omega_r})e^{-j\omega_r t}\right\}dt\le 0,
	\end{equation}
	which is equivalent to 
	\begin{equation}\label{eq:th5-3}
	\sum_{r=1}^{N} \hbox{Re}\left\{\lambda_rG({j\omega_r})\right\} 
	- \sum_{r=1}^{N-1} \hbox{Re}\left\{\lambda_rG({j\omega_r})H(j\omega_r)\right\}\le 0.
	\end{equation}
	Moreover, $H(j\infty)=0$ as $H$ is a strictly proper transfer function. Then we can write (\ref{eq:th5-3}) as
	\begin{equation}\label{eq:th5-4}
	\sum_{r=1}^{N} \hbox{Re}\left\{\lambda_rG({j\omega_r})\right\} 
	- \sum_{r=1}^{N} \hbox{Re}\left\{\lambda_rG({j\omega_r})H(j\omega_r)\right\}\le 0.
	\end{equation}
	Hence
	\begin{equation}\label{eq:th5-5}
	\sum_{r=1}^{N} \hbox{Re}\left\{\lambda_rM(j\omega_r)G({j\omega_r})\right\}\le 0,
	\end{equation}
	for all $M\in\mathcal{M}$, where $M=1-H$. Then the result is obtained by Corollary~\ref{cor:1}.
\end{IEEEproof}

\begin{IEEEbiography}[{\includegraphics[width=1in,height=1.25in,clip,keepaspectratio]{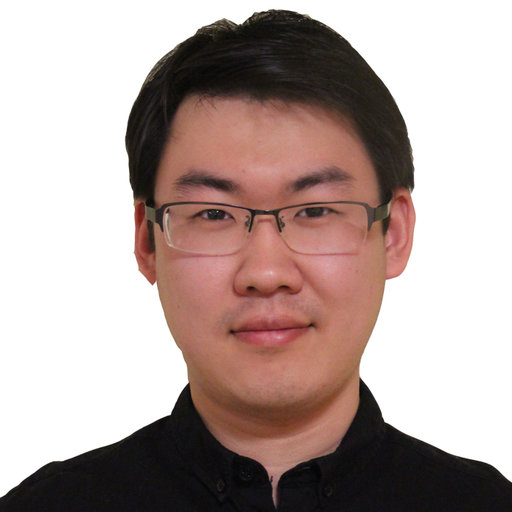}}]{Jingfan Zhang}
	received the BEng degree in electrical engineering and its automation from Xi'an Jiaotong-Liverpool University, Suzhou, China, in 2015, and MSc degree in advanced control and systems engineering from the University of Manchester, Manchester, U.K., in 2016. He is currently working toward the Ph.D. degree in the Department of Electrical and Electronic Engineering, University of Manchester, Manchester, U.K.
	His research interests include absolute stability and multiplier theory.
\end{IEEEbiography}

\begin{IEEEbiography}[{\includegraphics[width=1in,height=1.25in,clip,keepaspectratio]{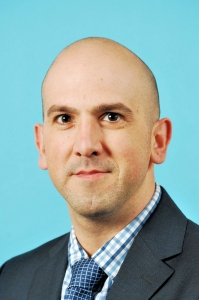}}]{Joaquin Carrasco}
	is a Senior Lecturer at the Control Systems Centre, Department of Electrical and Electronic Engineering, University of Manchester, UK. He was born in Abarán, Spain, in 1978. He received the B.Sc. degree in physics and the Ph.D. degree in control engineering from the University of Murcia, Murcia, Spain, in 2004 and 2009, respectively. From 2009 to 2010, he was with the Institute of Measurement and Automatic Control, Leibniz Universität Hannover, Hannover, Germany. From 2010 to 2011, he was a research associate at the Control Systems Centre, School of Electrical and Electronic Engineering, University of Manchester, UK. His current research interests include absolute stability, multiplier theory, and robotics applications. 
\end{IEEEbiography}

\begin{IEEEbiography}[{\includegraphics[width=1in,height=1.25in,clip,keepaspectratio]{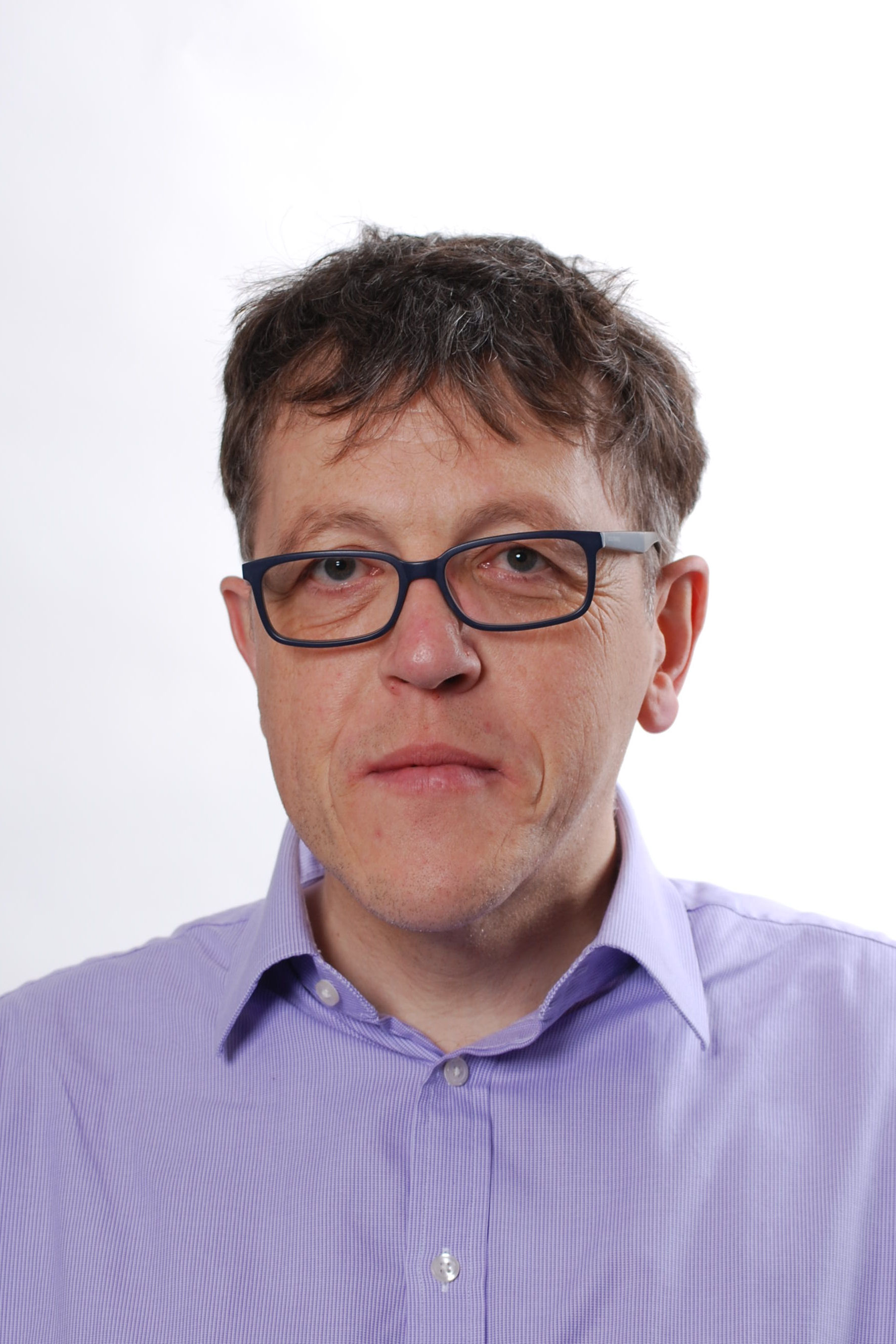}}]{William P. Heath}
	received an M.A. in mathematics from the University of Cambridge, U.K. and both an M.Sc. and Ph.D. in systems and control from the University of Manchester Institute of Science and Technology, U.K.
	He is Chair of Feedback and Control with the Control Systems Centre and Head of the Department of Electrical and Electronic Engineering, University of Manchester, U.K.  Prior to joining the University of Manchester, Professor Heath worked at Lucas Automotive and was a Research Academic at the University of Newcastle, Australia. His research interests include absolute stability, multiplier theory, constrained control, and system identification.
\end{IEEEbiography}

\end{algorithm}

\end{document}